\documentclass[12pt]{amsart}
\usepackage{amsmath, amsthm, amssymb}
\usepackage{tikz,tkz-graph}
\usetikzlibrary{shapes}
\usetikzlibrary{shapes.geometric}
\usepackage[top=1.25in, bottom=1.25in, left=1.0in, right=1.0in]{geometry}
\usepackage{hyperref}
\usepackage{color}
\usepackage{verbatim}
\usepackage{url}
\usepackage{subfig}
\newcounter{claim}
\makeatletter
\newtheorem*{rep@theorem}{\rep@title}
\newcommand{\newreptheorem}[2]{
\newenvironment{rep#1}[1]{
 \def\rep@title{#2 \ref{##1}}
 \begin{rep@theorem}}
 {\end{rep@theorem}}}
\makeatother

\theoremstyle{plain}
\newtheorem{thm}{Theorem}
\newreptheorem{thm}{Theorem}

\newreptheorem{prop}{Proposition}
\newtheorem{lem}[thm]{Lemma}
\newreptheorem{lem}{Lemma}

\newreptheorem{lemma}{Lemma}

\newreptheorem{conj}{Conjecture}

\newreptheorem{cor}{Corollary}
\newtheorem{cl}[claim]{Claim}

\theoremstyle{definition}

\theoremstyle{remark}

\newtheorem*{question}{Question}

\newcommand{\fancy}[1]{\mathcal{#1}}

\renewcommand{\Re}{\mathbb{R}}

\newcommand{\dist}{\textrm{dist}}
\newcommand{\Tiles}{\fancy{T}}
\newcommand{\IN}{\mathbb{N}}
\newcommand{\IR}{\mathbb{R}}
\newcommand{\Q}{\mathbb{Q}}
\newcommand{\ch}{\rm{ch}}
\newcommand{\set}[1]{\left\{ #1 \right\}}

\newcommand{\parens}[1]{\left( #1 \right)}

\def\naive{na\"{\i}ve}

\def\Naively{Na\"{\i}vely}

\def\aftermath{\par\vspace{-\belowdisplayskip}\vspace{-\parskip}\vspace{-\baselineskip}}

\tikzstyle{_RedDotStyle} = [shape = circle, minimum size = 10pt, draw, fill=black!30!red]
\tikzstyle{_GreenEdgeStyle} = [black!40!green, line width=3]
\tikzstyle{_BlueShadeStyle} = [shape = diamond, yscale=1.732, minimum size =
28pt, inner sep = 1.2pt, fill = blue!30]
\tikzstyle{_BlueShadeUpRightStyle} = [shape = regular polygon, regular polygon
sides=3, shape border rotate = 60, minimum size = 32.55pt, inner sep = 1.2pt, fill = blue!30]
\tikzstyle{_BlueShadeDownLeftStyle} = [shape = regular polygon, regular polygon
sides=3, shape border rotate = 240, minimum size = 32.55pt, inner sep = 1.2pt, fill = blue!30]
\tikzstyle{_ArrowStyle} = [->, ultra thick, >= triangle 45, blue!70]

\title{The Fractional Chromatic Number of the Plane}
\author{Daniel W. Cranston \and Landon Rabern}

\begin{document}
\begin{abstract}
The chromatic number of the plane is the chromatic number of the uncountably
infinite graph that has as its vertices the points of the plane and has an edge
between two points if their distance is 1.  This chromatic number is
denoted $\chi(\Re^2)$.  
The problem was introduced in 1950, and shortly thereafter 
it was proved that $4\le \chi(\Re^2)\le 7$.
These bounds are both easy to prove, but after more than 60 years they
are still the best known.  In this paper, we
investigate $\chi_f(\Re^2)$, the fractional chromatic number of the plane.
The previous best bounds (rounded to five decimal places) were $3.5556 \le
\chi_f(\Re^2)\le 4.3599$.  Here we improve the lower bound to $76/21\approx3.6190$.
\end{abstract}
\maketitle

\section{Introduction}

A proper coloring of the plane assigns to each of its points a color, such
that points at distance 1 get distinct colors.  The smallest number of colors that
allows such a coloring is the \emph{chromatic number of the plane}, denoted
$\chi(\Re^2)$.  This problem was introduced in 1950, by Edward Nelson, a student
at the University of Chicago.  
In the same year John Isbell, a fellow student, observed that $\chi(\Re^2)\le
7$.  This upper bound comes from a result of Hadwiger~\cite{Hadwiger45}, who showed that
it was possible to partition the plane into hexagons, each of radius slightly
less than 1, and color each hexagon with one of seven colors, such that hexagons
with the same color are distance more than 1 apart (see Figure~\ref{7coloring}).

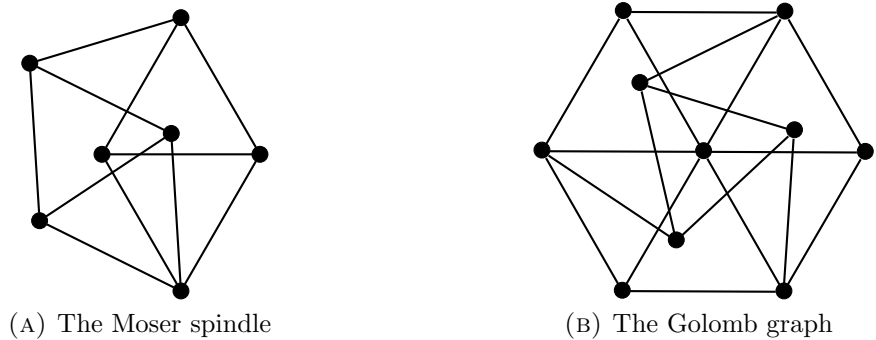
\begin{figure}
\subfloat[The Moser spindle]{\makebox[.45\textwidth]{
\begin{tikzpicture}[scale = 15]
\tikzstyle{VertexStyle}=[shape = circle, minimum size = 6pt, inner sep = 1.2pt,fill, draw]
\Vertex[x = 0.309031741140216, y = 0.872698767334361, L = \tiny {}]{v0}
\Vertex[x = 0.247631741140216, y = 0.854198767334361, L = \tiny {}]{v1}
\Vertex[x = 0.317631741140216, y = 0.975398767334361, L = \tiny {}]{v2}
\Vertex[x = 0.183631741140216, y = 0.934998767334361, L = \tiny {}]{v3}
\Vertex[x = 0.387631741140216, y = 0.854198767334361, L = \tiny {}]{v4}
\Vertex[x = 0.192331741140216, y = 0.795298767334361, L = \tiny {}]{v5}
\Vertex[x = 0.317631741140216, y = 0.732998767334361, L = \tiny {}]{v6}
\Edge[](v2)(v1)
\Edge[](v3)(v0)
\Edge[](v3)(v2)
\Edge[](v4)(v1)
\Edge[](v4)(v2)
\Edge[](v5)(v0)
\Edge[](v5)(v3)
\Edge[](v6)(v0)
\Edge[](v6)(v1)
\Edge[](v6)(v4)
\Edge[](v6)(v5)
\end{tikzpicture}}}
%
\subfloat[The Golomb graph]{\makebox[.45\textwidth]{
\begin{tikzpicture}[scale = 8,rotate=30]
\tikzstyle{VertexStyle}=[shape = circle, minimum size = 6pt, inner sep = 1.2pt,fill, draw]
\Vertex[x = 0.374613429155512, y = 0.82553315798011, L = \tiny {}]{v0}
\Vertex[x = 0.607716761976315, y = 0.958927166778096, L = \tiny {}]{v1}
\Vertex[x = 0.373727682881427, y = 0.556892033072417, L = \tiny {}]{v2}
\Vertex[x = 0.606743383611899, y = 0.690336276245073, L = \tiny {}]{v3}
\Vertex[x = 0.839856605170946, y = 0.823774121848011, L = \tiny {}]{v4}
\Vertex[x = 0.605785708356044, y = 0.421596564809465, L = \tiny {}]{v5}
\Vertex[x = 0.838890215317753, y = 0.555024109477077, L = \tiny {}]{v6}
\Vertex[x = 0.571808362602276, y = 0.841456096925729, L = \tiny {}]{v7}
\Vertex[x = 0.493244855598436, y = 0.584530428248699, L = \tiny {}]{v8}
\Vertex[x = 0.754976997291929, y = 0.644962349303816, L = \tiny {}]{v9}
\Edge[](v1)(v0)
\Edge[](v2)(v0)
\Edge[](v3)(v0)
\Edge[](v3)(v1)
\Edge[](v3)(v2)
\Edge[](v4)(v1)
\Edge[](v4)(v3)
\Edge[](v5)(v2)
\Edge[](v5)(v3)
\Edge[](v6)(v3)
\Edge[](v6)(v4)
\Edge[](v6)(v5)
\Edge[](v8)(v7)
\Edge[](v9)(v7)
\Edge[](v9)(v8)
\Edge[](v0)(v8)
\Edge[](v4)(v7)
\Edge[](v5)(v9)
\end{tikzpicture}}}
\caption{Two 4-chromatic unit distance graphs.\label{Moser-Golomb}}
\end{figure}

The lower bound $\chi(\Re^2)\ge 4$ comes from the observation of 
William and Leo Moser~\cite{Moser61} that the graph in Figure~\ref{Moser-Golomb}(A) is a unit
distance graph, i.e.,
it can be drawn in the plane with all edges of length 1. This graph is now known
as the Moser spindle; since it has chromatic number 4, the lower bound follows.
Around the same time~\cite[p.~19]{SoiferBook}, Solomon Golomb discovered the unit distance graph in
Figure~\ref{Moser-Golomb}(B), which also has chromatic number 4.
The problem first appeared in print in 1960\ in Martin Gardner's
\emph{Mathematical Games} column~\cite{Gardner60}.  

A seeming helpful result of de Bruijn and Erd\H{o}s~\cite{dBE51} implies that the chromatic
number of the plane is achieved by some finite subgraph (this proof does assume
the Axiom of Choice).   But unfortunately we have no reason to expect that such
a subgraph will have fewer than (say) a billion vertices.  By constructing a
$6$-coloring of nearly all of the plane, Pritikin \cite{pritikin1998} showed that if 
a $7$-chromatic unit distance graph does exist, then it has at least $6198$ vertices.
As for the lower bound, Erd\H{o}s wrote in 1985~\cite[p.~4]{erdos85} ``I am almost sure that
$\chi(\Re^2) > 4$.''

The history of the problem has many more interesting twists than we can recount
here, but Soifer records nearly all of them in his comprehensive and
entertaining \emph{The Mathematical Coloring Book}~\cite{SoiferBook}.
Although this problem has been widely popularized in the last half century, the
best known bounds remain $4\le \chi(\Re^2)\le 7$.

\input{7coloring}

The notion of fractional chromatic number was introduced in the early
1970s, with the goal of amassing more evidence in support of the Four Color
Conjecture, or possibly disproving it.  In a fractional coloring of a graph
$G$, we assign to each independent set in $G$ a nonnegative weight, such that each vertex
appears in independent sets with weights summing to at least 1.  The
\emph{fractional chromatic number}, $\chi_f(G)$, is the minimum sum of weights on the
independent sets that allows such a coloring.
This definition comes from solving the linear relaxation of the integer
programming formulation of chromatic number.

In 1992, Fisher and Ullman~\cite{FU92} first investigated the fractional chromatic number
of the plane.  They observed that the fractional chromatic number of the Moser
spindle is $3.5$, which gives a lower bound on $\chi_f(\Re^2)$.  They also gave
a coloring that proved the upper bound $\chi_f(\Re^2)\le 8\sqrt{3}/\pi \approx
4.4106$.  Using a similar approach, Hochberg and O'Donnell~\cite{HO93} improved
the upper bound to $4.3599$.  The construction they used was actually discovered much earlier by Croft \cite{croft1967}.
The lower bound was first improved by Shawna Mahan~\cite{Mahan95},
who found a unit distance graph with fractional chromatic number
$144/41\approx 3.5122$.  This bound was significantly improved by Fisher and
Ullman~\cite[p.~63--66]{SUbook}, who found a unit distance graph with
fractional chromatic number $\frac{32}9\approx 3.5556$.  The bounds $3.5556\le
\chi_f(\Re^2)\le 4.3599$ were the best known, until now.  In this paper, we
improve the lower bound to $76/21\approx3.6190$.

\section{A First Lower Bound}
\label{first-bound}
In this section we provide a unit distance graph with fractional
chromatic number greater than 3.6.  Our construction builds heavily on an
example of Fisher and Ullman, so we present that as well.  To begin, we
consider the fractional chromatic number of the two unit distance graphs we
have already seen, the Moser spindle and the Golomb graph.  The Moser spindle
has 7 vertices and independence number 2, which show that $\chi_f\ge \frac72=3.5$. 
To prove that this lower bound holds with equality, it suffices to
find 7 independent sets such that each vertex appears in two of them.  This
task is straightforward, once we put the bottom vertex into independent sets
with each of its nonneighbors.  The Golomb graph has 10 vertices and independence
number 3, which show that $\chi_f\ge \frac{10}3=3.\overline{3}$.  In
fact, this bound also holds with equality.  The matching upper bound comes
from finding 10 independent sets, such that each vertex appears in 3 of them. 
We leave this as an easy exercise.

The Moser spindle shows that $\chi_f(\Re^2)\ge 3.5$.  The intuition behind the
Fisher--Ullman construction is the following.  Consider a unit distance graph
that contains many copies of the Moser spindle, along with as many edges as
possible between these copies of the spindle.  These edges between
the copies of the spindle, as well as vertices that appear in more than one
copy, ensure that some copy of the spindle must be colored suboptimally.  This
will prove some lower bound greater than $3.5$.  The details forthwith. 

Recall that for any assignment of weights to the vertices,
the fractional chromatic number is bounded below by the total weight on the
vertices divided by the maximum total weight of any independent set.  (We
used this argument above to bound $\chi_f$ for the Moser spindle and Golomb
graph; there we implicitly gave each vertex weight 1.) Thus, for any weight
assignment, to bound $\chi_f$ from below, we need only bound from above the
maximum weight of any independent set.

\begin{figure}
\subfloat[The core vertices and weights from the Fisher--Ullman construction;
spindle weight 1 gives $\chi_f\ge \frac{32}{9} \approx 3.5555.$]{\makebox[.45\textwidth]{
\begin{tikzpicture}[scale = 10]
\tikzstyle{VertexStyle}=[shape = circle, minimum size = 6pt, inner sep = 1.2pt, draw]
\Vertex[x = 0.419762219286658, y = 0.9319, L = \tiny {$3$}]{v0}
\Vertex[x = 0.559762219286658, y = 0.9319, L = \tiny {$3$}]{v1}
\Vertex[x = 0.349762219286658, y = 0.8106, L = \tiny {$4$}]{v2}
\Vertex[x = 0.489762219286658, y = 0.8106, L = \tiny {$7$}]{v3}
\Vertex[x = 0.629762219286658, y = 0.8106, L = \tiny {$4$}]{v4}
\Vertex[x = 0.279762219286658, y = 0.6894, L = \tiny {$3$}]{v5}
\Vertex[x = 0.419762219286658, y = 0.6894, L = \tiny {$7$}]{v6}
\Vertex[x = 0.559762219286658, y = 0.6894, L = \tiny {$7$}]{v7}
\Vertex[x = 0.699762219286658, y = 0.6894, L = \tiny {$3$}]{v8}
\Vertex[x = 0.349762219286658, y = 0.5681, L = \tiny {$3$}]{v9}
\Vertex[x = 0.489762219286658, y = 0.5681, L = \tiny {$4$}]{v10}
\Vertex[x = 0.629762219286658, y = 0.5681, L = \tiny {$3$}]{v11}
\Edge[](v1)(v0)
\Edge[](v2)(v0)
\Edge[](v3)(v0)
\Edge[](v3)(v1)
\Edge[](v3)(v2)
\Edge[](v4)(v1)
\Edge[](v4)(v3)
\Edge[](v5)(v2)
\Edge[](v6)(v2)
\Edge[](v6)(v3)
\Edge[](v6)(v5)
\Edge[](v7)(v3)
\Edge[](v7)(v4)
\Edge[](v7)(v6)
\Edge[](v8)(v4)
\Edge[](v8)(v7)
\Edge[](v9)(v5)
\Edge[](v9)(v6)
\Edge[](v10)(v6)
\Edge[](v10)(v7)
\Edge[](v10)(v9)
\Edge[](v11)(v7)
\Edge[](v11)(v8)
\Edge[](v11)(v10)
\end{tikzpicture}}}
%
\subfloat[Vertices and weights for a bigger core;
spindle weight 1 gives $\chi_f\ge\frac{168}{47} \approx 3.5744.$]{\makebox[.45\textwidth]{
\begin{tikzpicture}[scale = 10]
\tikzstyle{VertexStyle}=[shape = circle, minimum size = 6pt, inner sep = 1.2pt, draw]
\Vertex[x = 0.571991341991342, y = 0.914744372294372, L = \tiny {$3$}]{v0}
\Vertex[x = 0.711991341991342, y = 0.914744372294372, L = \tiny {$3$}]{v1}
\Vertex[x = 0.501991341991342, y = 0.793544372294372, L = \tiny {$4$}]{v2}
\Vertex[x = 0.641991341991342, y = 0.793544372294372, L = \tiny {$7$}]{v3}
\Vertex[x = 0.781991341991342, y = 0.793544372294372, L = \tiny {$4$}]{v4}
\Vertex[x = 0.431991341991342, y = 0.672244372294372, L = \tiny {$4$}]{v5}
\Vertex[x = 0.571991341991342, y = 0.672244372294372, L = \tiny {$8$}]{v6}
\Vertex[x = 0.711991341991342, y = 0.672244372294372, L = \tiny {$8$}]{v7}
\Vertex[x = 0.851991341991342, y = 0.672244372294372, L = \tiny {$4$}]{v8}
\Vertex[x = 0.361991341991342, y = 0.551044372294372, L = \tiny {$3$}]{v9}
\Vertex[x = 0.501991341991342, y = 0.551044372294372, L = \tiny {$7$}]{v10}
\Vertex[x = 0.641991341991342, y = 0.551044372294372, L = \tiny {$8$}]{v11}
\Vertex[x = 0.781991341991342, y = 0.551044372294372, L = \tiny {$7$}]{v12}
\Vertex[x = 0.921991341991342, y = 0.551044372294372, L = \tiny {$3$}]{v13}
\Vertex[x = 0.431991341991342, y = 0.429844372294372, L = \tiny {$3$}]{v14}
\Vertex[x = 0.571991341991342, y = 0.429844372294372, L = \tiny {$4$}]{v15}
\Vertex[x = 0.711991341991342, y = 0.429844372294372, L = \tiny {$4$}]{v16}
\Vertex[x = 0.851991341991342, y = 0.429844372294372, L = \tiny {$3$}]{v17}
\Edge[](v1)(v0)
\Edge[](v2)(v0)
\Edge[](v3)(v0)
\Edge[](v3)(v1)
\Edge[](v3)(v2)
\Edge[](v4)(v1)
\Edge[](v4)(v3)
\Edge[](v5)(v2)
\Edge[](v6)(v2)
\Edge[](v6)(v3)
\Edge[](v6)(v5)
\Edge[](v7)(v3)
\Edge[](v7)(v4)
\Edge[](v7)(v6)
\Edge[](v8)(v4)
\Edge[](v8)(v7)
\Edge[](v9)(v5)
\Edge[](v10)(v5)
\Edge[](v10)(v6)
\Edge[](v10)(v9)
\Edge[](v11)(v6)
\Edge[](v11)(v7)
\Edge[](v11)(v10)
\Edge[](v12)(v7)
\Edge[](v12)(v8)
\Edge[](v12)(v11)
\Edge[](v13)(v8)
\Edge[](v13)(v12)
\Edge[](v14)(v9)
\Edge[](v14)(v10)
\Edge[](v15)(v10)
\Edge[](v15)(v11)
\Edge[](v15)(v14)
\Edge[](v16)(v11)
\Edge[](v16)(v12)
\Edge[](v16)(v15)
\Edge[](v17)(v12)
\Edge[](v17)(v13)
\Edge[](v17)(v16)
\end{tikzpicture}}}
\caption{The two smallest cores that we consider.\label{small-cores}}
\end{figure}
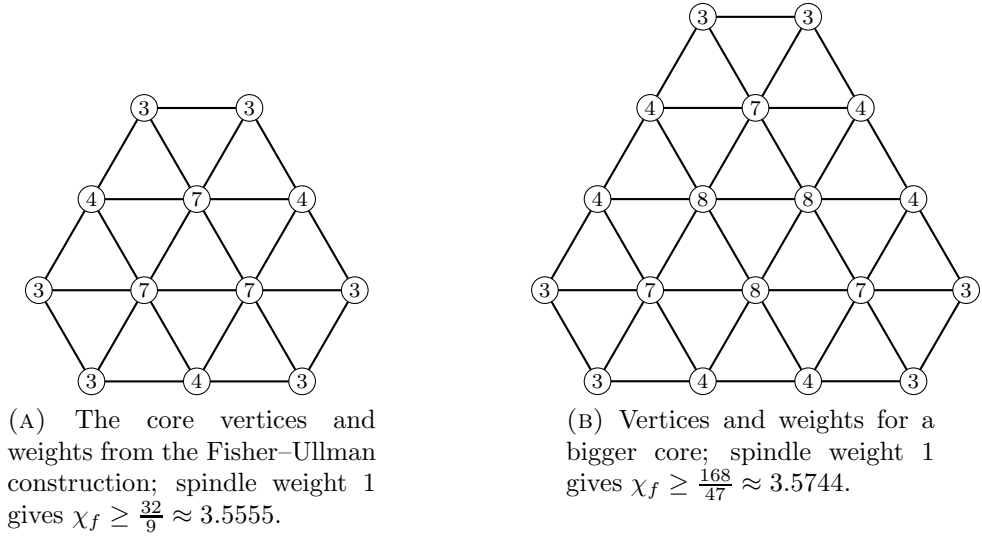

We construct the Fisher--Ullman graph in two stages.  We begin with the subset
of the triangular lattice shown in Figure~\ref{small-cores}(A), called the
\emph{core}; for now ignore the weights, which we will get to shortly.  
A \emph{diamond} is the subgraph induced by two vertices at distance $\sqrt{3}$
and their two common neighbors.  For each of the 5 vertical diamonds, we attach
a copy of the Moser spindle; in each case, we identify the four vertices of the
diamond with the four vertices of the vertical diamond in the spindle, and we
add three new vertices.  These new vertices are \emph{spindle vertices}.
The two core vertices that are each adjacent to at least one of these spindle
vertices are \emph{incident to the spindle}.  Nothing is special about the
vertical direction in the core, so we also attach spindles to the five diamonds
pointing down and to the left, as well as the five pointing down and to the
right; before attaching these spindles, we rotate them 120 degrees clockwise
and 120 degrees counterclockwise, respectively.  Finally, we add all edges
between pairs of vertices at distance 1.  

The graph that results has 3-fold
rotational symmetry.  Each spindle adds 3 new vertices, so the 15 spindles add
a total of 45 vertices.  With the 12 core vertices, this makes a total of
57 vertices.  The graph has 24 edges among core vertices and $6(15)$ more edges within
spindles; for each of three directions, it has $21$ edges between vertices in
different spindles that are oriented in the same direction.  Finally, it has
another 21 edges among pairs that seem to ``accidentally'' be at distance 1. 
(Understanding these last 21 edges is inessential, since they are not used in
proving the lower bound.)

Now we assign weights to the vertices.  The weights on the core vertices are
shown in Figure~\ref{small-cores}(A).  To each spindle vertex, we assign weight
1.  We have 45 spindle vertices, and the weights on the core sum to 51, so the
total weight is 96.  Thus, to prove a lower bound of $\frac{32}9$, it suffices
to show that every independent set has weight at most $27$.  This requires a
short case analysis, and Scheinerman and Ullman \cite[p.~64--65]{SUbook} gives
most of the details.

Following this proof is fairly easy, but where do the weights come from?
They come from solving a linear program (hereafter LP).  Specifically, each weight is a
variable, each independent set has weight at most 1, and the sum of the weights
is to be maximized. (To simplify our presentation above, we multiplied all
weights by 27, but that does not affect the lower bound.)
So what about larger cores?

\begin{figure}
\subfloat[The core vertices and weights for a bigger core;
spindle weight 2 gives $\chi_f\ge\frac{491}{137} \approx 3.5839$.]{\makebox[.45\textwidth]{
\begin{tikzpicture}[scale = 10]
\tikzstyle{VertexStyle}=[shape = circle, minimum size = 6pt, inner sep = 1.2pt, draw]
\Vertex[x = 0.721355389541089, y = 0.897550373532551, L = \tiny {$5$}]{v0}
\Vertex[x = 0.861355389541089, y = 0.897550373532551, L = \tiny {$5$}]{v1}
\Vertex[x = 0.651355389541089, y = 0.776350373532551, L = \tiny {$6$}]{v2}
\Vertex[x = 0.791355389541089, y = 0.776350373532551, L = \tiny {$12$}]{v3}
\Vertex[x = 0.931355389541089, y = 0.776350373532551, L = \tiny {$6$}]{v4}
\Vertex[x = 0.581355389541089, y = 0.655050373532551, L = \tiny {$7$}]{v5}
\Vertex[x = 0.721355389541089, y = 0.655050373532551, L = \tiny {$16$}]{v6}
\Vertex[x = 0.861355389541089, y = 0.655050373532551, L = \tiny {$16$}]{v7}
\Vertex[x = 1.00135538954109, y = 0.655050373532551, L = \tiny {$7$}]{v8}
\Vertex[x = 0.511355389541089, y = 0.533850373532551, L = \tiny {$6$}]{v9}
\Vertex[x = 0.651355389541089, y = 0.533850373532551, L = \tiny {$16$}]{v10}
\Vertex[x = 0.791355389541089, y = 0.533850373532551, L = \tiny {$20$}]{v11}
\Vertex[x = 0.931355389541089, y = 0.533850373532551, L = \tiny {$16$}]{v12}
\Vertex[x = 1.07135538954109, y = 0.533850373532551, L = \tiny {$6$}]{v13}
\Vertex[x = 0.441355389541089, y = 0.412550373532551, L = \tiny {$5$}]{v14}
\Vertex[x = 0.581355389541089, y = 0.412550373532551, L = \tiny {$12$}]{v15}
\Vertex[x = 0.721355389541089, y = 0.412550373532551, L = \tiny {$16$}]{v16}
\Vertex[x = 0.861355389541089, y = 0.412550373532551, L = \tiny {$16$}]{v17}
\Vertex[x = 1.00135538954109, y = 0.412550373532551, L = \tiny {$12$}]{v18}
\Vertex[x = 1.14135538954109, y = 0.412550373532551, L = \tiny {$5$}]{v19}
\Vertex[x = 0.511355389541089, y = 0.291350373532551, L = \tiny {$5$}]{v20}
\Vertex[x = 0.651355389541089, y = 0.291350373532551, L = \tiny {$6$}]{v21}
\Vertex[x = 0.791355389541089, y = 0.291350373532551, L = \tiny {$7$}]{v22}
\Vertex[x = 0.931355389541089, y = 0.291350373532551, L = \tiny {$6$}]{v23}
\Vertex[x = 1.07135538954109, y = 0.291350373532551, L = \tiny {$5$}]{v24}
\Edge[](v1)(v0)
\Edge[](v2)(v0)
\Edge[](v3)(v0)
\Edge[](v3)(v1)
\Edge[](v3)(v2)
\Edge[](v4)(v1)
\Edge[](v4)(v3)
\Edge[](v5)(v2)
\Edge[](v6)(v2)
\Edge[](v6)(v3)
\Edge[](v6)(v5)
\Edge[](v7)(v3)
\Edge[](v7)(v4)
\Edge[](v7)(v6)
\Edge[](v8)(v4)
\Edge[](v8)(v7)
\Edge[](v9)(v5)
\Edge[](v10)(v5)
\Edge[](v10)(v6)
\Edge[](v10)(v9)
\Edge[](v11)(v6)
\Edge[](v11)(v7)
\Edge[](v11)(v10)
\Edge[](v12)(v7)
\Edge[](v12)(v8)
\Edge[](v12)(v11)
\Edge[](v13)(v8)
\Edge[](v13)(v12)
\Edge[](v14)(v9)
\Edge[](v15)(v9)
\Edge[](v15)(v10)
\Edge[](v15)(v14)
\Edge[](v16)(v10)
\Edge[](v16)(v11)
\Edge[](v16)(v15)
\Edge[](v17)(v11)
\Edge[](v17)(v12)
\Edge[](v17)(v16)
\Edge[](v18)(v12)
\Edge[](v18)(v13)
\Edge[](v18)(v17)
\Edge[](v19)(v13)
\Edge[](v19)(v18)
\Edge[](v20)(v14)
\Edge[](v20)(v15)
\Edge[](v21)(v15)
\Edge[](v21)(v16)
\Edge[](v21)(v20)
\Edge[](v22)(v16)
\Edge[](v22)(v17)
\Edge[](v22)(v21)
\Edge[](v23)(v17)
\Edge[](v23)(v18)
\Edge[](v23)(v22)
\Edge[](v24)(v18)
\Edge[](v24)(v19)
\Edge[](v24)(v23)
\end{tikzpicture}}}
%
\subfloat[The core vertices and weights for a bigger core; spindle weight 3
gives $\chi_f\ge\frac{1732}{481} \approx 3.6008$.  This LP had over $25,000$ constraints.]{\makebox[.45\textwidth]{
\begin{tikzpicture}[scale = 15]
\tikzstyle{VertexStyle}=[shape = circle, minimum size = 6pt, inner sep = 1.2pt, draw]
\Vertex[x = 0.491985974754558, y = 0.928522230014025, L = \tiny {$6$}]{v0}
\Vertex[x = 0.547985974754558, y = 0.928522230014025, L = \tiny {$6$}]{v1}
\Vertex[x = 0.463985974754558, y = 0.880022230014025, L = \tiny {$11$}]{v2}
\Vertex[x = 0.519985974754558, y = 0.880022230014025, L = \tiny {$21$}]{v3}
\Vertex[x = 0.575985974754558, y = 0.880022230014025, L = \tiny {$11$}]{v4}
\Vertex[x = 0.435985974754558, y = 0.831522230014025, L = \tiny {$9$}]{v5}
\Vertex[x = 0.491985974754558, y = 0.831522230014025, L = \tiny {$26$}]{v6}
\Vertex[x = 0.547985974754558, y = 0.831522230014025, L = \tiny {$26$}]{v7}
\Vertex[x = 0.603985974754558, y = 0.831522230014025, L = \tiny {$9$}]{v8}
\Vertex[x = 0.407985974754558, y = 0.783022230014025, L = \tiny {$9$}]{v9}
\Vertex[x = 0.463985974754558, y = 0.783022230014025, L = \tiny {$19$}]{v10}
\Vertex[x = 0.519985974754558, y = 0.783022230014025, L = \tiny {$21$}]{v11}
\Vertex[x = 0.575985974754558, y = 0.783022230014025, L = \tiny {$19$}]{v12}
\Vertex[x = 0.631985974754558, y = 0.783022230014025, L = \tiny {$9$}]{v13}
\Vertex[x = 0.379985974754558, y = 0.734522230014025, L = \tiny {$9$}]{v14}
\Vertex[x = 0.435985974754558, y = 0.734522230014025, L = \tiny {$18$}]{v15}
\Vertex[x = 0.491985974754558, y = 0.734522230014025, L = \tiny {$18$}]{v16}
\Vertex[x = 0.547985974754558, y = 0.734522230014025, L = \tiny {$18$}]{v17}
\Vertex[x = 0.603985974754558, y = 0.734522230014025, L = \tiny {$18$}]{v18}
\Vertex[x = 0.659985974754558, y = 0.734522230014025, L = \tiny {$9$}]{v19}
\Vertex[x = 0.351985974754558, y = 0.686022230014025, L = \tiny {$9$}]{v20}
\Vertex[x = 0.407985974754558, y = 0.686022230014025, L = \tiny {$19$}]{v21}
\Vertex[x = 0.463985974754558, y = 0.686022230014025, L = \tiny {$18$}]{v22}
\Vertex[x = 0.519985974754558, y = 0.686022230014025, L = \tiny {$19$}]{v23}
\Vertex[x = 0.575985974754558, y = 0.686022230014025, L = \tiny {$18$}]{v24}
\Vertex[x = 0.631985974754558, y = 0.686022230014025, L = \tiny {$19$}]{v25}
\Vertex[x = 0.687985974754558, y = 0.686022230014025, L = \tiny {$9$}]{v26}
\Vertex[x = 0.323985974754558, y = 0.637622230014025, L = \tiny {$11$}]{v27}
\Vertex[x = 0.379985974754558, y = 0.637622230014025, L = \tiny {$26$}]{v28}
\Vertex[x = 0.435985974754558, y = 0.637622230014025, L = \tiny {$21$}]{v29}
\Vertex[x = 0.491985974754558, y = 0.637622230014025, L = \tiny {$18$}]{v30}
\Vertex[x = 0.547985974754558, y = 0.637622230014025, L = \tiny {$18$}]{v31}
\Vertex[x = 0.603985974754558, y = 0.637622230014025, L = \tiny {$21$}]{v32}
\Vertex[x = 0.659985974754558, y = 0.637622230014025, L = \tiny {$26$}]{v33}
\Vertex[x = 0.715985974754558, y = 0.637622230014025, L = \tiny {$11$}]{v34}
\Vertex[x = 0.295985974754558, y = 0.589122230014025, L = \tiny {$6$}]{v35}
\Vertex[x = 0.351985974754558, y = 0.589122230014025, L = \tiny {$21$}]{v36}
\Vertex[x = 0.407985974754558, y = 0.589122230014025, L = \tiny {$26$}]{v37}
\Vertex[x = 0.463985974754558, y = 0.589122230014025, L = \tiny {$19$}]{v38}
\Vertex[x = 0.519985974754558, y = 0.589122230014025, L = \tiny {$18$}]{v39}
\Vertex[x = 0.575985974754558, y = 0.589122230014025, L = \tiny {$19$}]{v40}
\Vertex[x = 0.631985974754558, y = 0.589122230014025, L = \tiny {$26$}]{v41}
\Vertex[x = 0.687985974754558, y = 0.589122230014025, L = \tiny {$21$}]{v42}
\Vertex[x = 0.743985974754558, y = 0.589122230014025, L = \tiny {$6$}]{v43}
\Vertex[x = 0.323985974754558, y = 0.540622230014025, L = \tiny {$6$}]{v44}
\Vertex[x = 0.379985974754558, y = 0.540622230014025, L = \tiny {$11$}]{v45}
\Vertex[x = 0.435985974754558, y = 0.540622230014025, L = \tiny {$9$}]{v46}
\Vertex[x = 0.491985974754558, y = 0.540622230014025, L = \tiny {$9$}]{v47}
\Vertex[x = 0.547985974754558, y = 0.540622230014025, L = \tiny {$9$}]{v48}
\Vertex[x = 0.603985974754558, y = 0.540622230014025, L = \tiny {$9$}]{v49}
\Vertex[x = 0.659985974754558, y = 0.540622230014025, L = \tiny {$11$}]{v50}
\Vertex[x = 0.715985974754558, y = 0.540622230014025, L = \tiny {$6$}]{v51}
\Edge[](v1)(v0)
\Edge[](v2)(v0)
\Edge[](v3)(v0)
\Edge[](v3)(v1)
\Edge[](v3)(v2)
\Edge[](v4)(v1)
\Edge[](v4)(v3)
\Edge[](v5)(v2)
\Edge[](v6)(v2)
\Edge[](v6)(v3)
\Edge[](v6)(v5)
\Edge[](v7)(v3)
\Edge[](v7)(v4)
\Edge[](v7)(v6)
\Edge[](v8)(v4)
\Edge[](v8)(v7)
\Edge[](v9)(v5)
\Edge[](v10)(v5)
\Edge[](v10)(v6)
\Edge[](v10)(v9)
\Edge[](v11)(v6)
\Edge[](v11)(v7)
\Edge[](v11)(v10)
\Edge[](v12)(v7)
\Edge[](v12)(v8)
\Edge[](v12)(v11)
\Edge[](v13)(v8)
\Edge[](v13)(v12)
\Edge[](v14)(v9)
\Edge[](v15)(v9)
\Edge[](v15)(v10)
\Edge[](v15)(v14)
\Edge[](v16)(v10)
\Edge[](v16)(v11)
\Edge[](v16)(v15)
\Edge[](v17)(v11)
\Edge[](v17)(v12)
\Edge[](v17)(v16)
\Edge[](v18)(v12)
\Edge[](v18)(v13)
\Edge[](v18)(v17)
\Edge[](v19)(v13)
\Edge[](v19)(v18)
\Edge[](v20)(v14)
\Edge[](v21)(v14)
\Edge[](v21)(v15)
\Edge[](v21)(v20)
\Edge[](v22)(v15)
\Edge[](v22)(v16)
\Edge[](v22)(v21)
\Edge[](v23)(v16)
\Edge[](v23)(v17)
\Edge[](v23)(v22)
\Edge[](v24)(v17)
\Edge[](v24)(v18)
\Edge[](v24)(v23)
\Edge[](v25)(v18)
\Edge[](v25)(v19)
\Edge[](v25)(v24)
\Edge[](v26)(v19)
\Edge[](v26)(v25)
\Edge[](v27)(v20)
\Edge[](v28)(v20)
\Edge[](v28)(v21)
\Edge[](v28)(v27)
\Edge[](v29)(v21)
\Edge[](v29)(v22)
\Edge[](v29)(v28)
\Edge[](v30)(v22)
\Edge[](v30)(v23)
\Edge[](v30)(v29)
\Edge[](v31)(v23)
\Edge[](v31)(v24)
\Edge[](v31)(v30)
\Edge[](v32)(v24)
\Edge[](v32)(v25)
\Edge[](v32)(v31)
\Edge[](v33)(v25)
\Edge[](v33)(v26)
\Edge[](v33)(v32)
\Edge[](v34)(v26)
\Edge[](v34)(v33)
\Edge[](v35)(v27)
\Edge[](v36)(v27)
\Edge[](v36)(v28)
\Edge[](v36)(v35)
\Edge[](v37)(v28)
\Edge[](v37)(v29)
\Edge[](v37)(v36)
\Edge[](v38)(v29)
\Edge[](v38)(v30)
\Edge[](v38)(v37)
\Edge[](v39)(v30)
\Edge[](v39)(v31)
\Edge[](v39)(v38)
\Edge[](v40)(v31)
\Edge[](v40)(v32)
\Edge[](v40)(v39)
\Edge[](v41)(v32)
\Edge[](v41)(v33)
\Edge[](v41)(v40)
\Edge[](v42)(v33)
\Edge[](v42)(v34)
\Edge[](v42)(v41)
\Edge[](v43)(v34)
\Edge[](v43)(v42)
\Edge[](v44)(v35)
\Edge[](v44)(v36)
\Edge[](v45)(v36)
\Edge[](v45)(v37)
\Edge[](v45)(v44)
\Edge[](v46)(v37)
\Edge[](v46)(v38)
\Edge[](v46)(v45)
\Edge[](v47)(v38)
\Edge[](v47)(v39)
\Edge[](v47)(v46)
\Edge[](v48)(v39)
\Edge[](v48)(v40)
\Edge[](v48)(v47)
\Edge[](v49)(v40)
\Edge[](v49)(v41)
\Edge[](v49)(v48)
\Edge[](v50)(v41)
\Edge[](v50)(v42)
\Edge[](v50)(v49)
\Edge[](v51)(v42)
\Edge[](v51)(v43)
\Edge[](v51)(v50)
\end{tikzpicture}}}
\caption{Two larger cores. \label{big-cores}}
\end{figure}
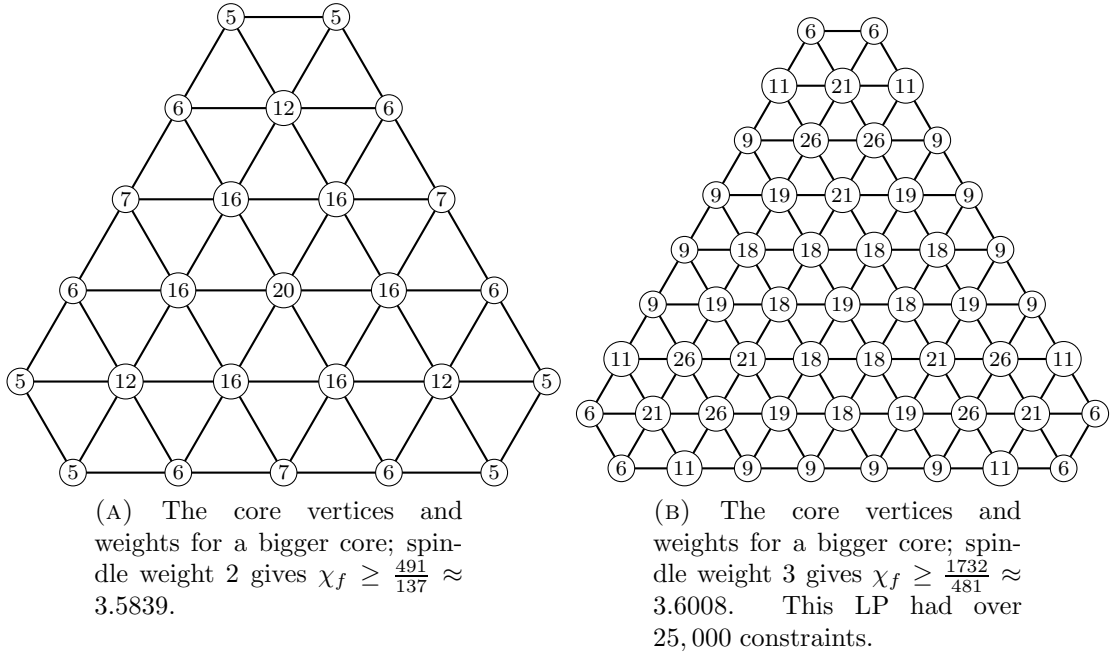

We can easily generalize the Fisher--Ullman construction to start from a larger
core (as suggested in~\cite[p.75]{SUbook}), and we will do exactly this.  But
first, it is helpful to comment on the obvious symmetry of the weights in
Figure~\ref{small-cores}(A).  Suppose we are given an optimal assignment of
weights, i.e., an optimal solution to the LP in the previous paragraph. 
Each automorphism of the graph yields another optimal assignment of weights.
Further, the average of all these weight assignments is again optimal.
Thus, we may assume that the same weight is given to all vertices in
each orbit when the vertices are acted on by the automorphism group.  This
observation~\cite[p.7]{LevinThesis} dramatically reduces the size of our LP,
thus making it tractable.

Figures~\ref{small-cores}(B) and \ref{big-cores}(A) show the results from
solving the corresponding LPs for the next few sizes of cores.  
In each case, the lower bound on $\chi_f$ improves, but more slowly.  
We set for ourselves the goal of proving a lower bound of 3.6, and using the
weights in Figure~\ref{big-cores}(B) we achieved it, just barely.
In fact, these weights prove $\chi_f\ge 3.6008$.
We would have liked to consider still larger cores, but the size of the derived
LP was growing exponentially (along with the number of maximal independent
sets), and the LP for the core in Figure~\ref{big-cores}(B) involved already
more than 25,000 constraints.  For this lower bound of 3.6008, we offer no
proof.  However, the interested reader can download our code for generating the LP from: \url{https://github.com/landon/WebGraphs/blob/master/Analysis/SpindleAnalyzer.cs}.
To do much better, and certainly to give a human-checkable proof, we needed a
new approach.

\section{An Improved Lower Bound}
\label{improved-bound}
\subsection{A Lower Bound by Discharging}
\label{first-discharging}

In the previous section, we proved that $\chi_f(\Re^2)>3.6$.  To do so, we chose
a unit distance graph and gave its vertices weights summing to more than 3.6,
such that the weights on every independent set summed to at most 1.
However, this proof has three drawbacks.  First, it is not practically
human-checkable, since that graph has more than 25,000 maximal independent sets.
So to verify that each independent set has weight at most 1, we have a lot of
cases to check.  Second, we have no reason to believe that this bound of 3.6008
is actually near the fractional chromatic number of the plane.  It is simply the
best bound we could prove before the number of maximal independent sets became
unmanageable.  Finally, and perphaps most disturbingly, the proof offers no real
insight.  We have simply found some weights ``that work''.

In this section, we prove a stronger lower bound.  In the process, we address
all of these concerns, as well.  We take a similar approach to our previous
proof, but with a few key differences.  As before, we start with a subset of
the triangular lattice, called the core.  Now everywhere that we possibly can, we add
on spindles, much like in the previous proof.  The first main difference is
that we don't worry much about optimizing the weights.  We give the same weight
to every vertex in the core, and the same weight to every spindle vertex.  (We
will optimize the ratio of these two weights, but that problem is much
easier.)  

The second main difference---really the key that allows the proof to
work---is that when bounding
the weight of an independent set, we don't consider the set all at once.
Rather, for an arbitrary maximal independent set $I$, we partition the graph into
subgraphs, and bound the fraction of the weight on each subgraph that is in 
$I$.  Now the fraction of weight on the whole graph that is $I$ is no more than
the maximum fraction on any of these subgraphs.  
By partitioning into subgraphs of bounded size, we avoid the combinatorial
explosion we faced in the previous section, when each maximal independent set
generated its own constraint in the LP.
An important decision is how to choose these subgraphs, so that we have
relatively few cases, but we also get a good bound on the fraction of the total
weight in $I$.  We return to this question later.

To illustrate our approach, we start with a short proof that $\chi_f(\Re^2)\ge 3.5$.
Obviously this bound is weak (recall that the Moser spindle has fractional
chromatic number 3.5), but its proof elucidates many of the key features of our
method.
\bigskip

We define a graph $G_d$ and its weights as follows.
To begin, specify an arbitrary vertex $v$ in the triangular lattice.
For our core, take all of the triangular lattice induced by vertices that are
distance at most $d$ from $v$.  For our spindles, add every possible spindle in
each of the three directions we used in the previous section.  This is our
graph $G_d$.  Let $C_d$ denote the subgraph of $G_d$ induced by its core
vertices.  For our weight function, we assign weight 12 to every core vertex
and weight 1 to every spindle vertex.  Now, given an arbitrary independent set
$I$, we will discharge all of its weight to core vertices, so that each core
vertex has weight at most 6.  We use the following two discharging rules.

\begin{enumerate}
\item[(R1)] Each core vertex in $I$ gives weight 1 to each of its
neighbors in the core.
\item[(R2)] Each spindle vertex in $I$ splits it weight
equally between the core vertices incident to its spindle that are
\emph{not} in $I$.
\end{enumerate}

Note that if a spindle vertex is in $I$, then at least one of the core vertices
incident to the spindle is not in $I$; hence $(R2)$ does, indeed, move all
charge from spindle vertices to core vertices.
Now clearly all the weight in $I$ ends on vertices in the core.  We must verify
that each core vertex finishes with charge at most 6.
In additition to core vertices in $I$, we have four possibilites for a core vertex
not in $I$; it can have 3, 2, 1, or 0 core neighbors in $I$.  We now check these
five cases.

Note that each core vertex has at most 6 incident spindles (exactly 6 if the
core vertex is not too close to the ``outside'' of the core).  So a \naive\ upper
bound on the final charge at each core vertex $v$ is 9, since $v$ receives 1
from each of at most 3 core neighbors, and it receives at most 1 from each of 6
spindles.  However, this bound can be improved.  Suppose that $v$ is a core
vertex not
in $I$, and let $u$ be some core neighbor of $v$ that is in $I$.  Now $u$ has
two neighbors, say $w_1$ and $w_2$, that are each distance $\sqrt{3}$
from $v$.  Note that $v$ shares one spindle each with $w_1$ and $w_2$.
The key observation is that $u\in I$ implies $w_1\notin I$ and $w_2\notin I$. 
Thus, $v$ receives weight at most $1/2$ from each of the spindles it shares
with $w_1$ and $w_2$.  Repeatedly applying this insight leads to the following
upper bounds on the final charge at each core vertex.

\textbf{Core vertex in $I$:} $12-6(1)+0 = 6$.

\textbf{Core vertex with three neighbors in $I$:} $0 + 3(1) + 6(1/2) = 6$.  

\textbf{Core vertex with two neighbors in $I$:} $0  + 2(1) + 4(1/2) + 2(1) = 6$.

\textbf{Core vertex with one nbr in $I$:} $0 + 1(1) + 2(1/2) + 4(1) = 6$.

\textbf{Core vertex with zero neighbors in $I$:} $0 + 0(1) + 0(1/2) + 6(1) = 6$.

So we have shown that each core vertex finishes with weight at most $6$.
To compute a lower bound on $\chi_f$, we divide the total weight on
$G_d$ by (an upper bound on) the total weight of any independent set.
To simplify the computations, we neglect the effect of vertices near the outside
of the core, which have fewer than six incident spindles (this choice can be
justified, since the number of ``interior'' core vertices is asymptotically
greater than the number of ``boundary'' core vertices).  Let $M$ denote the
number of core vertices.  Each core vertex is incident to six spindles and each
spindle is incident to two core vertices; hence, the number of spindles is
$M(3-o(1))$.  Since each spindle has 3 spindle vertices, the number of spindle
vertices is $M(9-o(1))$.  Thus, the total weight on the graph is $M(12+9(1)-o(1))$.
The final total weight on the core is at most $6M$.  This proves a lower
bound of $\frac{21}6 =3.5$.

It is reasonable to ask why we give the same weight to every core vertex, since this
differs from our choice in the previous section.  A simple answer is that ``it
works'', but this is unenlightening.  A somewhat better answer is that  ``by
weighting each vertex identically, we dramatically reduce the number of cases we
must consider.''  This is true, but still misses the real point.  To understand
more fully, we return to the weights used in Section 2.  There we weighted two
vertices identically whenever the graph had an automorphism mapping one to the
other.  In that context, our automorphisms consisted of reflections and
rotations, and compositions of the two.  

In the present context, since we
consider $G_d$ as $d$ grows without bound, we are essentially considering the
fractional chromatic number of an infinite graph.  (Our choice to consider the
graph in ever larger pieces, rather than all at once, is mainly a concession
with the goal of simplifying the averaging argument.)  If we are indeed
considering an infinite graph, then we have an additional type of automorphism:
translations.  Now it is clear that we can map any core vertex to any other core
vertex; so the choice to weight them equally is quite natural.

Once we understand the proof of this lower bound, it's natural to ask how we
can strengthen it.  One obvious choice is to change the weight we give to each
core vertex.  To make this approach work, we need to consider the core vertices
in larger groups, rather than in isolation, as in the previous proof.  Our idea
is to partition the subgraph induced by the core vertices into smaller
subgraphs, which we call \emph{tiles}.  As our tiles grow bigger, we
gain more information about the neighborhoods of their vertices, which
facilitates a better bound on their average final weight.  However, as the tiles
grow bigger, they also grow more numerous, and require more case analysis.  As a
compromise, we choose the tiles to be (essentially) as small as possible,
subject to each corner of each tile being a vertex of $I$, the independent set.

A priori, the distance could be large between a core vertex in $I$ and its nearest
core vertex in $I$.  However, we are interested only in choices of $I$ that
potentially have maximum weight.  By requiring that each core vertex gets
weight greater than its number of incident spindles, we can make the following
simplification.  We only need to consider choices of $I$ that are maximal
independent subsets of the core.  Suppose instead that for some core
vertex $v$ some independent set $I$ contains neither $v$ nor any of its core
neighbors.  We form $I'$ from $I$ by adding $v$ and removing all spindle
neighbors of $v$ in $I$.  Since the weight on $v$ is greater than its number of
incident spindles, $I'$ weighs more than $I$.

In the following lemma, we formalize this idea of choosing the tiles to be ``as
small as possible, subject to each corner of each tile being a vertex of $I$,
the independent set.''  We also show that this process results in only 8
distinct tiles, up to rotation and reflection.

\subsection{Tiling Result}
\label{tiling}

\input{figures}

\begin{lem}
Let $I$ denote a maximal independent subset in $C_d$.  There exists a set
$\Tiles$ of 8 finite \emph{tiles} (shown in Figure~\ref{tiles}), independent
of $d$ and $I$, such that $C_d$ can be tiled with tiles from $\Tiles$ where
each corner of each tile is a vertex of $I$ and no vertex of $I$ lies in
the interior of any tile. In this tiling, each face of $C_d$ is covered by
exactly one or two tiles.  (We do allow tiles to extend past the boundary of
$G_d$, though this allowance could be removed by adding more tiles
to $\Tiles$.)  

\end{lem}
\begin{proof}
For an example of such a tiling, see Figure~\ref{tiling-pic}.

Before we begin the formal proof, we note that every maximal subset $I$ of $C_d$
must contain at least $\frac17$ of its vertices.  Since $I$ is maximal, for every
vertex $v$, the set $I$ contains either $v$ or some neighbor of $v$.  So to
prove the lower bound $\frac17$, we observe that there exists a set $A$ of
$\frac17$ of all the vertices such that for every pair $u,w\in S$, we have
$\dist_{C_d}(u,w) \ge 3$.  Thus, for each vertex $u$ in $A$, the set $I$
contains some element of $N[u]$; further, these sets are disjoint, i.e., for
any pair $u,w\in A$, we have $N[u]\cap N[w] = \emptyset$.  To see that such a
set $A$ exists, we view the infinite triangular lattice as the planar dual of
the 7-face-coloring shown in Figure~\ref{7coloring}.  This gives a 7-coloring
of the triangular lattice, and hence a 7-coloring of $C_d$; we choose as $A$
the vertices of a largest color class.

This reasoning actually proves that $I$ must contain at least
roughly $\frac17$ of the vertices in any region of $C_d$.  Intuitively, for each
face $f$ of $C_d$, some nearby vertices are in $I$, so it will be possible to
cover $f$ using a small tile, with all of its vertices nearby.  Now we make this
intuition rigorous.

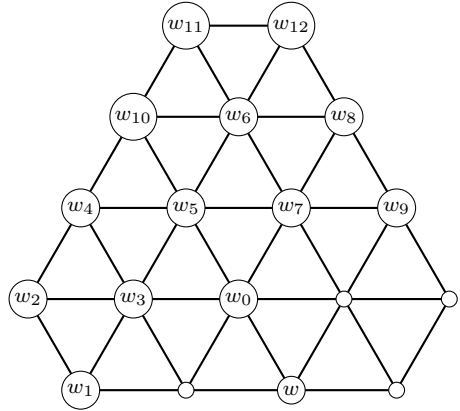
\begin{figure}[hbt]
\begin{tikzpicture}[scale = 5]
\tikzstyle{VertexStyle}=[shape = circle, minimum size = 6pt, inner sep = 1.2pt, draw]

\begin{scope}[xshift=1.55cm,yshift=-0.24248711306cm]
\Vertex[x = 1.21302593659942, y = 0.7822163824, L = \tiny {$w_{11}$}]{w28}
\Vertex[x = 1.49302593659942, y = 0.7822163824, L = \tiny {$w_{12}$}]{w29}
\Vertex[x = 1.07302593659942, y = 0.539729269346951, L = \tiny {$w_{10}$}]{w18}
\Vertex[x = 1.35302593659942, y = 0.539729269346951, L = \tiny {$w_6$}]{w2}
\Vertex[x = 1.63302593659942, y = 0.539729269346951, L = \tiny {$w_8$}]{w13}
\Vertex[x = 0.933025936599424, y = 0.297242156287308, L = \tiny {$w_4$}]{w12}
\Vertex[x = 1.21302593659942, y = 0.297242156287308, L = \tiny {$w_5$}]{w1}
\Vertex[x = 1.49302593659942, y = 0.297242156287308, L = \tiny {$w_7$}]{w3}
\Vertex[x = 1.77302593659942, y = 0.297242156287308, L = \tiny {$w_9$}]{w4}
\Vertex[x = 0.793025936599424, y = 0.0547550432276658, L = \tiny {$w_2$}]{w17}
\Vertex[x = 1.07302593659942, y = 0.0547550432276658, L = \tiny {$w_3$}]{w11}
\Vertex[x = 1.35302593659942, y = 0.0547550432276658, L = \tiny {$w_0$}]{w0}
\Vertex[x = 1.63302593659942, y = 0.0547550432276658, L = \tiny {}]{w5}
\Vertex[x = 1.91302593659942, y = 0.0547550432276658, L = \tiny {}]{w14}
\Vertex[x = 0.933025936599424, y = -0.187732069831977, L = \tiny {$w_1$}]{w10}
\Vertex[x = 1.21302593659942, y = -0.187732069831977, L = \tiny {}]{w9}
\Vertex[x = 1.49302593659942, y = -0.187732069831977, L = \tiny {$w$}]{w7}
\Vertex[x = 1.77302593659942, y = -0.187732069831977, L = \tiny {}]{w6}
\Edge[](w0)(w1)
\Edge[](w0)(w3)
\Edge[](w0)(w5)
\Edge[](w0)(w7)
\Edge[](w0)(w9)
\Edge[](w0)(w11)
\Edge[](w1)(w2)
\Edge[](w1)(w3)
\Edge[](w1)(w11)
\Edge[](w1)(w12)
\Edge[](w1)(w18)
\Edge[](w2)(w3)
\Edge[](w2)(w13)
\Edge[](w2)(w18)
\Edge[](w3)(w4)
\Edge[](w3)(w5)
\Edge[](w3)(w13)
\Edge[](w4)(w5)
\Edge[](w4)(w13)
\Edge[](w4)(w14)
\Edge[](w5)(w6)
\Edge[](w5)(w7)
\Edge[](w5)(w14)
\Edge[](w6)(w7)
\Edge[](w6)(w14)
\Edge[](w7)(w9)
\Edge[](w9)(w10)
\Edge[](w9)(w11)
\Edge[](w10)(w11)
\Edge[](w10)(w17)
\Edge[](w11)(w12)
\Edge[](w11)(w17)
\Edge[](w12)(w17)
\Edge[](w12)(w18)
\Edge[](w28)(w29)
\Edge[](w28)(w2)
\Edge[](w28)(w18)
\Edge[](w29)(w2)
\Edge[](w29)(w13)
\end{scope}

\end{tikzpicture}
\caption{Names of key vertices for verifying that the tiles have only 8 shapes.
\label{names}}
\end{figure}

First we describe the process for partitioning $C_d$ into tiles; shortly, we will 
analyze it to show that we form only 8 tiles (up to rotation and reflection).  
To partition $C_d$ into pieces, for each vertex pair $u,w\in I$, we add edge
$uw$ if and only if $u$ and $w$ have euclidean distance less than 3.  Whenever
two edges cross, we delete them both.  This process clearly constructs a plane
graph.  In what follows, we show that the faces of this graph, which will
become our tiles, have at most 8 distinct shapes and sizes.
Let $T$ be some arbitrary face of the plane graph constructed above.  We first 
consider the case where $T$ contains some edge of $C_d$ in its interior (not
boundary) that is incident to one of its corners; we will later show that if
this is not the case, then $T$ must be T2.
By symmetry, we may assume that this corner is $w$, as shown in
Figure~\ref{names}.  
By rotational symmetry around $w$, we may assume that $T$
contains edge $ww_0$ in its interior.

As noted above, $I$ contains $w_3$ or one of its neighbors; similarly, $I$
contains $w_7$ or one of its neighbors.  Since $ww_0$ is in the interior of
$T$, we may assume that $w_5\notin I$.  Since $w\in I$, none of its neighbors
are in $I$.  Thus, $I$ contains some vertex in each of vertex sets $A$ and $B$, where
$A=\set{w_1,w_2,w_3,w_4}$ and $B=\set{w_6,w_7,w_8,w_9}$.
\Naively, this gives 16 cases to consider. 
The number is actually smaller due to the symmetry across
the line segment $\overline{ww_5}$, but we rarely make use of this symmetry,  
to keep the case organization clear. 
We consider these cases now.  

Suppose that $w_3\in I$.
If $w_7\in I$, then we have T1 (throughout the proof, we refer to each tile by
its caption in Figure~\ref{tiles}).  If $w_8\in I$, then we have T4.
If exactly one of $w_6$ and $w_9$ is in $I$, then we have T3.  If both of
$w_6$ and $w_9$ are in $I$, then we have T6; note that this is the unique
instance where initially two edges of our constructed graph crossed, so we
deleted them.

Suppose that $w_4\in I$.  If $w_6\in I$, then we have T4.  If $w_7\in I$,
then we have T3.  If $w_8\in I$, then we have T5.  So, suppose that
$w_9\in I$ and $w_6\notin I$.  Recall that $I$ must contain $w_6$ or one of its
neighbors.  So either $w_{11}\in I$ or $w_{12}\in I$.  In the former case, we
have T7; in the latter, we have T8.

Suppose that $w_2\in I$.  If $w_6\in I$, then we have T5.  If $w_7\in I$,
then we have T4.  If $w_8\in I$, then we must have $w_{10}\in I$, since
$I$ contains neither $w_5$ nor any of its other neighbors.  Now we have
T8.  So suppose that $w_9\in I$ and $w_6\notin I$.  Again, we must have
$w_{10}\in I$, since $I$ contains neither $w_5$ nor any of its other neighbors.
Now we have T7.

Suppose that $w_1\in I$ (and $w_4\notin I$).  By symmetry across
$\overline{ww_5}$, we may assume that $w_9\in I$ and $w_6\notin I$.  Again,
this implies $w_{10}\in I$.  Once again, we have T8.

Now we consider the case that $T$ contains no edge of $C_d$ in its interior.
Let $w$ be a corner of $T$.  We observe, as follows, that each boundary segment
of $T$ must have length at most 2.  Suppose, to the contrary, that
$T$ has a boundary segment that is longer.  By symmetry, we may assume that it
is $\overline{ww_2}$.  By symmetry, we may also assume that $T$ lies above this
segment (the case where $T$ lies below it is essentially the same).  We will
show that $T$ contains the edge $w_2w_3$ in its interior.  Since $I$ is
maximal, it contains $w_5$ or one of its neighbors; specifically, $I$ contains
$w_5$, $w_6$, $w_7$, or $w_{10}$.  In the first three cases, $T$ is a triangle
and clearly contains $w_2w_3$ in its interior.  So assume that $w_{10}\in I$
and $w_7\notin I$.  Now either $w_8\in I$ or $w_9\in I$.  So $T$ is either
T8 or T7, respectively; in each case $T$ contains edge $w_2w_3$ in
its interior.  Hence, we conclude that each boundary segment of $T$ has
length at most 2.

Suppose now that $T$ has a boundary segment of length $\sqrt{3}$.
By symmetry, say this is segment $\overline{ww_3}$ and that $T$ lies above this
segment.  Since $I$ is maximal, it contains either $w_7$ or (at least) one of
its neighbors.  In each case, $T$ contains edge $ww_0$ in its interior.
Thus, we conclude that $T$ has no boundary segments of length
$\sqrt{3}$.  Hence, every boundary segment of $T$ must have length 2.
Further, the interior angle formed by boundary segments at a corner must be
$\pi/3$.  So we conclude that $T$ must be T2.
\end{proof}

\input{tiling}

The previous result shows that we can tile $C_d$ with only 8 tiles.  However,
each of these tiles may be rotated or reflected in numerous ways.  In order to
discharge the weight from spindle vertices to core vertices, we will view each
of these tiles together with the spindles incident to its vertices.  If we keep
with our earlier strategy of attaching spindles at each vertex $v$ in only 3
directions with $v$ as their bottom vertex (for a total of 6 incident
spindles, including those with $v$ as their top vertex), we will break some of the
symmetries we required in our proof that 8 tiles suffices.  In particular, the
roles of vertices incident to a spindle at its top and bottom are inherently
asymmetrical, since the bottom is adjacent to two spindle vertices, while the
top is adjacent to only one.  

One obvious approach for handling this difficulty is to consider 
multiple cases for each tile, depending on its orientation; but this solution is
inelegant.  So we prefer instead the following approach.  We now attach spindles
at each vertex in 6 directions (for a total of 12 incident spindles) as follows.  
Each pair of vertices, say $u$ and $v$, that shared a spindle before now share
two spindles; for one spindle $u$ is the top vertex and $v$ is the bottom, and
for the other spindle the roles are reversed.  For a pair of vertices $u$ and
$v$ that share a spindle, the positions of the spindle vertices in the second
spindle come from reflecting
the vertices in the first spindle across the perpendicular bisector of
line segment $\overline{uv}$.  Thus, if the spindle with $u$ as its base is rotated
$\theta$ radians clockwise (around $u$) past $\overline{uv}$ (for an
appropriate $\theta$), then the spindle with $v$ as its base is rotated
$\theta$ radians \emph{counterclockwise} (around $v$) past $\overline{uv}$.
Hence, each vertex $v$ will be the bottom vertex for 6 spindles. Three of these
six spindles will each be oriented $2\pi/3$ radians clockwise of the previous one;
the other three will be oriented similarly (relative to each other).
We call this graph $G'_d$.

In Section~\ref{first-discharging}, we used discharging to prove a lower bound
of 3.5.  To rephrase that proof in this setting with twice as many spindles, we
simply give each core vertex weight 12 and each spindle vertex weight $\frac12$
(rather than 1, as before).  A moment's reflection will show that the rest of
that proof goes through essentially unchanged.  

An attentive reader will perhaps wonder whether all of these supposedly distinct 
spindle vertices do indeed fall in different locations.  The answer is yes,
although it turns out not to matter.  First, we should mention that the angle
of rotation for each spindle, called $\theta$ above, is $\cos^{-1}(\frac56)$.  
In less than a page of computations, we can show that the spindle vertices
really do have distinct locations.  But we don't need to.

Recall that our proof that $\chi_f(\Re^2)\ge \frac{76}{21}$ consists of two
parts.  In one part, we prove lower bounds on the fractional chromatic number
for a sequences of graphs, and show that these lower bounds converge to
$\frac{76}{21}$.  In the other part, we show that each graph in the sequence is
a unit distance graph.  We have constructed our graph sequence $G'_d$, and in
what follows, we will prove the desired lower bounds.  So we need only show
that each $G'_d$ is a unit distance graph, by describing some
embedding.  

We take as our embedding of the core vertices the obvious one,
from the triangular lattice.  If it happens that two (or more) spindle vertices
coincide, then we assign to this ``combined'' vertex at their common location the
\emph{sum} of the weights we had intended to assign to them individually.  Each
of the supposedly distinct spindle vertices now has more neighbors than we
originally claimed, but none has fewer.  If this combined spindle vertex $v$
does appear in the independent set $I$, then, for each spindle $S$ containing
$v$, we discharge to the core vertices incident $S$ the portion of $v$'s weight
that was due to it appearing in $S$.  Thus, the lower bound remains valid.

\subsection{Discharging Result}

In this section, we continue with the graphs $G'_d$ of the previous section.  We
show that as $d$ grows, the fractional chromatic number of $G'_d$ is bounded
below by a sequence converging to $\frac{76}{21}\approx 3.619047$. In what
follows, when we write $o(1)$, we mean as $d$ goes to infinity.

\begin{thm}
The fractional chromatic number of the plane is at least $\frac{76}{21}$, i.e.,
$\chi_f(\Re^2)\ge \frac{76}{21}$.
\end{thm}

\begin{proof}
At each core vertex we attach spindles in 6 directions, as described above.
Each core vertex gets weight $\frac{31}5$, and each spindle vertex gets weight
$\frac12$.  Let $I$ be an arbitrary independent set.
We will redistribute the weight in $I$ so that the core vertices ends with average
weight at most $\frac{21}5$ and each spindle vertex ends with weight at most 0.
To help bound the average final weight of the core vertices, we compute the
average weight of the core vertices in each tile, and show that the average
weight for each tile is at most $\frac{21}5$.  As before, we let $M$ denote the
number of core vertices.  The total weight on the core is $\frac{31}5M$, and
each spindle vertex has weight $\frac12$.  Since the number of spindle vertices
is $M(18-o(1))$, the total weight on $G'_d$ is $M(\frac{31}5+\frac12(18)-o(1))$.
Since the core vertices end with average weight at most $\frac{21}5$, the total
weight in $I$ is at most $\frac{21}5M$.
This proves that $\chi_f(\Re^2)\ge
(\frac{31}5+\frac12(18)-o(1))/(\frac{21}5) =
\frac{76}{21}\approx 3.619047$.

Now we give the details of how to redistribute the weight.
Each vertex of $C_d$ has a target weight of (at most) $\frac{21}5$.
The \emph{target} weight for a tile $T$ is $\frac{21}5(i_T+\frac12b_T+0c_T)$, where
$i_T$, $b_T$, and $c_T$ denote the number of core vertices (respectively) in the interior
of $T$, on the boundary (but not corners) of $T$, and on the corners of $T$.
Each corner of $T$ is a vertex of $I$, and its final weight will be computed by
itself.  Each interior vertex of $T$  has all of its target weight assigned to
$T$; similarly, each boundary (but not corner) vertex of $T$ has exactly one
half of its weight assigned to $T$.  If $T$ has more weight than its target,
then this difference is its \emph{excess}; otherwise, the difference is a
\emph{deficit}. 
Our goal is to show that each tile finishes with excess at most 0.

Each diamond has a spindle attached in two directions; we call these the
\emph{up} spindle and the \emph{down} spindle.  If a diamond has two vertices in $I$,
then its spindles are \emph{trivial}.  A non-trivial spindle that has no spindle vertex
in $I$ is \emph{missing}.  We redistribute the charge using three discharging phases.

\subsubsection*{Phase 1}
\begin{enumerate}
\item[(R1)] Each core vertex in $I$ gives weight $\frac13$ to each core neighbor.

\item[(R2)] Each non-trivial spindle splits weight $\frac12$ equally
among the core vertices not in $I$ that are incident to the spindle.  If one of
the two incident core vertices is in $I$, then the spindle sends weight $\frac12$ to the
other; otherwise it sends weight $\frac14$ to each.

\item[(R3)] Each core vertex $v$ not in $I$ gives its weight to the tile
containing it.  If $v$ is on the border of two tiles, then $v$ splits
its weight between them (but not quite equally).  This is a little subtle,
but crucial.  For such a $v$, it has two neighbors in $I$; $v$ does split the
weight from its neighbors in $I$ equally among its two tiles.  Now $v$ also
has 12 incident spindles; so $v$ sends any weight it got from 6 of the
spindles to one tile and the weight from the other spindles to the
other tile.  Split the 12 spindles into two sets of 6 each, spaced 120
degrees apart (the up spindle and down spindle for each pair of vertices go into
the same set).  Now $v$ gives each tile the weight it got from the set
of spindles that includes the spindle pointing away from the border
into that tile.
\end{enumerate}

\subsubsection*{Phase 2}
\begin{enumerate}
\item[(R4)] 
Each copy of T3 takes $\frac3{10}$ from the tile on its long side.
\item[(R5)] Each copy of T1 splits its deficit equally among all adjacent copies of
T6.  In other words, it takes an equal amount of weight from each adjacent
copy of T6, so that it ends with excess 0.
\end{enumerate}

\subsubsection*{Phase 3}
\begin{enumerate}
\item[(R6)]
Each copy of T2 gives $\frac14$ to each missing spindle in
its 5-spindle block (Figure~\ref{spindle-blocks}(A)) in each of three directions.
Each copy of T6 gives $\frac18$ to each missing spindle in
its 6-spindle block (Figure~\ref{spindle-blocks}(B)) in each of four directions.
\end{enumerate}
\input{spindle-blocks}

\bigskip

We now check that each vertex in $I$ and each tile ends with average weight at most
$\frac{21}5$.  For vertices in $I$, this is immediate since each ends with 
weight $\frac{31}5 - 6(\frac13)=\frac{21}5$.  We also check that each spindle
finishes with weight at most 0.
\bigskip

\begin{cl}
\label{cl1}
After Phase 1 the tiles have excess at most:\\
T1: -1/5
T2: 7/10 
T3: -3/10
T4: 3/5
T5: 2/5
T6: 2/5
T7: 0
T8: 0.
\end{cl}
\begin{proof}
The main observation needed for this proof is essentially the same one we needed
for the discharging proof in Section~\ref{first-discharging}.
If $v$ is a core
vertex not in $I$, then for each core neighbor $u$ of $v$ in $I$ there exist four
spindles incident to $v$ that send $v$ weight at most $\frac14$ each.
Specifically, $u$ has two neighbors, say $w_1$ and $w_2$, that are each
distance $\sqrt{3}$ from $v$.  Further, $v$ shares two spindles with each $w_i$.
Since $u\in I$, we know that $w_1,w_2\notin I$.  Thus, each spindle shared
between $v$ and each $w_i$ splits its weight equally between $v$ and $w_i$.
Hence, $v$ receives weight at most $\frac14$ from each such spindle.

{T1:} $3(\frac13) +12(\frac14) = 4$.  The target weight is $\frac{21}5$, so the excess is $4-\frac{21}5=-\frac{1}5$.

{T2:} 
Since T2 has 3-fold rotational symmetry, we compute the weight received from one boundary
vertex, then multiply by 3.  This is
$3(\frac13+4(\frac14)+2(\frac12)) = 7$.
The target weight is
$\frac32(\frac{21}5)=\frac{63}{10}$,
so the excess is
$\frac{7}{10}$.

{T3:}
We may assume that the tile $\widehat{T}$ bordering this copy $T$ of T3 along its longest side is
not another copy of T3 with each of its edges parallel to an edge of $T$.  If it
were, then these two copies of T3 would be merged into a single copy of T6.
This assumption enables us to improve the bound on the total weight on $T$,
which is $(\frac13+6(\frac14))+(2(\frac13)+10(\frac14)+2(\frac12)) = 6$.  The
target weight is $\frac32(\frac{21}{5}) = \frac{63}{10}$, so the excess is
$-\frac{3}{10}$.

{T4:}
$(2(\frac13)+12(\frac14))+(\frac13+4(\frac14)+8(\frac12))=9$.  The target weight is
$2(\frac{21}5)=\frac{42}5$, so the excess is $\frac{3}5$.

{T5:}
Since T5 has 3-fold rotational symmetry, we compute the weight received from one
interior vertex, then multiply by 3. This is
$3(\frac13+8(\frac14)+4(\frac12))= 13$.  The target weight is $3(\frac{21}5)=\frac{63}5$,
so the excess is $\frac{2}5$.

{T6:}
Since T6 has 2-fold rotational symmetry, we compute the weight received from
vertices $u$ and $v$, then multiply by 2.  This is
$2((\frac13+6(\frac14))+(2(\frac13)+8(\frac14)+4(\frac12)))=13$.  The target weight is
$3(\frac{21}5)=\frac{63}5$, so the excess is $\frac{2}5$.

{T7:}
Since T7 has 2-fold rotational symmetry, we compute
the weight received from $u$, $v$, and $w$, then multiply by 2.  This is
$2((\frac13+6(\frac14))+(\frac13+6(\frac14)+6(\frac12))+(\frac13+10(\frac14)+2(\frac12)))=21$.  
The target weight is $5(\frac{21}5)=21$, so the excess is 0.

{T8:}
Since T8 has reflectional symmetry, we compute the weight received from $w$ and
$z$ plus twice the weight received from $u$ and $x$.
This is
$((\frac13+8(\frac14)+4(\frac12))+(\frac13+12(\frac14)))+2((\frac13+6(\frac14))+(\frac13+6(\frac14)+6(\frac12)))=21$.
The target weight is $5(\frac{21}5)=21$, so the excess is 0.
\end{proof}

\begin{cl}
\label{cl2}
After Phase~2 the tiles have excess at most:

T1: 0
T2: 7/10
T3: 0
T4: 0
T5: 0
T6: 2/5
T7: 0
T8: 0.
\end{cl}
\begin{proof}
Since Phase 2 only affects the weight for copies of T1, T3, and T6, the bounds
for the other Figures from Claim~1 remain valid.
By (R5), each copy of T1 ends with excess 0.  Note that T3 has only
a single long side.  Since each copy ended Phase 1 with excess at most $-\frac{3}{10}$,
it ends Phase 2 with excess at most 0.  

Consider a copy $T$ of T4, and the
tile $\widehat{T}$ that borders it along one of its two long sides.  If
$\widehat{T}$ is a copy
of T3, then T4 gives away $\frac3{10}$.  So assume instead that $\widehat{T}$ is a copy of
T4, T5, T7, or T8.
In this case, $T$ saves at least $\frac12$ over the bound we computed in
Claim~\ref{cl1}.  Specifically, an additional vertex at distance $\sqrt{3}$ from $v$ is
not in $I$.  This means that $v$ receives only one half of the weight from each
spindle with these two vertices as its top and bottom (rather than all of it).
So $T$ saves $\frac{3}{10}$ if $\widehat{T}$ is a copy of T3, and saves $\frac12$
otherwise.  This savings applies to both of the long edges of T4, so T4 saves at
least $2(\frac3{10})=\frac35$ over the bound given in Claim~\ref{cl1}.  Thus,
each copy of T4 finishes with excess at most 0.
A similar analysis holds for T5.  Now T5 has three long edges, and saves at
least $\frac3{10}$ along each of them, so finishes with excess at most
$\frac25-3(\frac3{10})=-\frac{1}2$.  
\end{proof}

\bigskip

\input{claim3pics}

\begin{cl}
\label{cl3}
Each spindle ends with nonnegative weight.
\end{cl}
\begin{proof}
If a spindle is trivial, then it begins and ends with weight 0.
If a spindle is non-trivial and not missing, then it begins with weight
$\frac12$.  By (R2), it ends with weight 0.
So we need only consider missing spindles.  By (R2), each such spindle finishes
Phases 1 and 2 with weight $-\frac{1}{2}$.
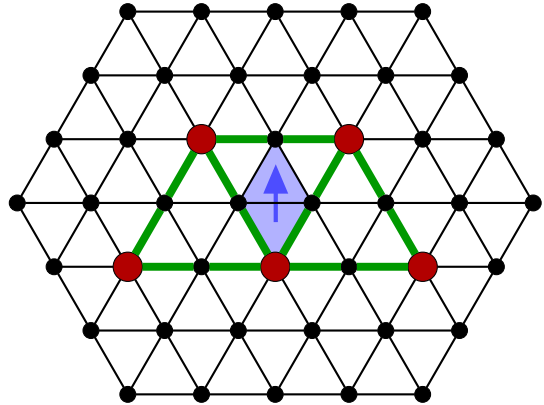
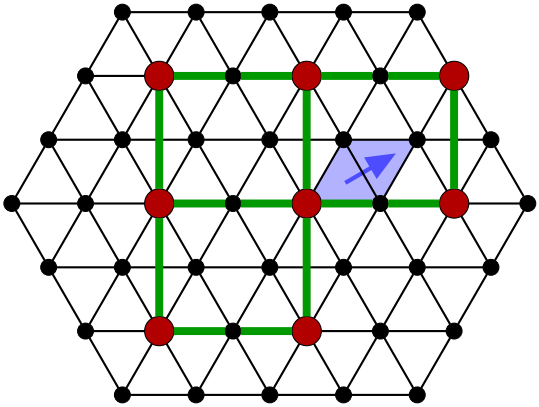
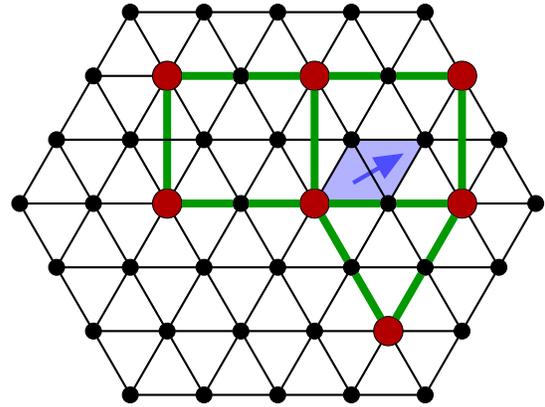
So it suffices to show that in Phase 3 each missing spindle receives weight at
most $\frac12$.  Recall that a missing spindle receives weight $\frac14$ (from a
copy of T2) for each 5-spindle block containing it and weight $\frac18$ (from a 
copy of T6) for each 6-spindle block containing it.  Thus, we must bound the
number of 5-spindle and 6-spindle blocks containing each missing spindle.

First, if a missing spindle lies in two 5-spindle blocks, then it lies in no
other 5-spindle or 6-spindle blocks; see Figure~\ref{cl3pics}(A).
So it receives weight at most $\frac14+\frac14= \frac12$.  
Further, each missing spindle lies in at most three 6-spindle blocks.
Such a spindle lies in no 5-spindle blocks; see Figure~\ref{cl3pics}(B).
So it receives weight at most $\frac18+\frac18+\frac18\le \frac12$. 
Finally, if a spindle lies in both a 5-spindle block and a 6-spindle block,
then it lies in at most one 5-spindle block and at most two 6-spindle blocks;
see Figure~\ref{cl3pics}(C).
So it receives weight at most $\frac18+\frac18+\frac14= \frac12$.
\end{proof}

\begin{cl}
In each 5-spindle block of a copy of T2, at least one spindle is missing.
\label{cl4}
\end{cl}
\begin{proof}
We number the (spindle) vertices of a spindle as 1, 2, 3, where vertices 1 and 2
are distance one from the bottom vertex, and vertex 3 is distance $\sqrt{3}$.
We number the spindles in a 5-spindle block as $S_1,\ldots,S_5$, from left to
right.  Suppose that none of $S_1,\ldots, S_5$ is missing.  Since the bottom
vertices of $S_1$ and $S_5$ are in $I$, each of them must also have vertex 3 in
$I$.  In each of $S_2, S_3, S_4$ vertex 3 is distance 1 from vertex 3
of either $S_1$ or $S_5$.  Since both of these vertices are in $I$, none of
$S_2, S_3, S_4$ can have vertex 3 in $I$.  However, among $S_2, S_3, S_4$, the
three copies of vertex 1 are pairwise adjacent; similarly for the three copies
of vertex 2.  Hence, one of $S_2, S_3, S_4$ must be missing.
\end{proof}

\begin{cl}
\label{cl5}
For each copy of T6, the following are true.
\begin{enumerate}
\item[(a)]
For each 6-spindle block,
either (i) at least one spindle is missing or (ii) at least two spindles are trivial.

\item[(b)]
If a copy of T6 has 2 trivial spindles on one side in the same direction, then that side
is bordered by a copy of T1 that is adjacent to at most two copies of T6.  

\item[(c)]
If a copy of T6 has 4 trivial spindles on one side, then that side
is bordered by a copy of T1 that is adjacent to only one copy of T6.

\end{enumerate}
\end{cl}
\begin{proof}
We prove each part in turn.

\begin{enumerate}
\item[(a)] The proof is similar to that of Claim~\ref{cl4}.
We number the spindles as $S_1,\ldots, S_6$, from left to right (with the
topmost spindle as $S_5$ and the bottom rightmost as $S_6$).  We number the
spindle vertices of each spindle 1, 2, 3, as in Claim~\ref{cl5}.
The only possible trivial spindles are $S_5$ and $S_6$; so suppose that at least
one of them is neither trivial nor missing.  Such a spindle must have vertex 3
in $I$.  Whether this is $S_5$ or $S_6$, it forbids $S_4$ from having vertex 3
in $I$.  Similarly, if $S_1$ is not missing, it has vertex 3 in $I$, which
forbids both of spindles $S_2$ and $S_3$ from having vertex 3 in $I$.  Now
spindles $S_2, S_3, S_4$ all have vertex 3 forbidden from $I$.  However, their
copies of vertex 1 (resp. vertex 2) are pairwise adjacent.  Hence, one of $S_2,
S_3, S_4$ is missing.
\item[(b)]
Recall from (a) that the only possible trivial spindles in the 6-spindle block
are $S_5$ and $S_6$.  If $S_6$ is trivial, then the copy of T6 is bordered on
its short side by a copy of T1.  Further, if $S_5$ is also trivial, then that
copy of T1 shares one of its sides with another copy of T1; hence, it can be
bordered by at most two copies of T6.
\item[(c)]
Each short side of a copy of T6 is crossed by two of its 6-spindle blocks. 
Recall that each 6-spindle block has only two possible trivial spindles.  Thus,
if a copy of T6 has 4 trivial spindles on one side, then both of its 6-spindle blocks
crossing that side must have both possible trivial spindles.  Now we apply (b)
twice.  So the copy of T6 is bordered along this side by a copy of T1, and that
copy of T1 is bordered along its two other sides by other copies of T1.  Hence,
the copy of T1 is adjacent to only one copy of T6.
\end{enumerate}
\aftermath
\vspace{.05in}~
\end{proof}

\begin{cl}
\label{cl6}
Each tile ends with excess at most 0.
\end{cl}
\begin{proof}
By Claim~2, we only need to check copies of T2 and T6.
Further, at the end of Phase 2, each copy of T2 has excess at most $\frac7{10}$
and each copy of T6 has excess at most $\frac25$.
By Claim~\ref{cl4}, each copy of
T2 gives weight $\frac14$ to a missing spindle in each of three
directions.  Thus, its final excess is at most
$\frac7{10}-3(\frac14)=-\frac{1}{20}$.

If a copy of T6 has a missing spindle in all 4 directions, then it gives away
to its missing spindles at least $4(\frac18) = \frac12 \ge \frac25$.  If T6 has
a missing spindle in 3 directions, then it gives away to its missing spindles 
at least $3(\frac18) = \frac38$.  By Claim~\ref{cl5}(a), in the direction where
there is no missing spindle, T6 has 2 trivial spindles on one side in the same
direction.  By Claim~\ref{cl5}(b), the copy of T6 is bordered by a copy of T1
that has at most 2 adjacent copies of T6.  Hence by (R5), the copy of T6 gives
at least $\frac12(\frac15) = \frac1{10}$ to the copy of T1.  Thus the copy of
T6 finishes with excess at most $\frac25-\frac38 -\frac1{10}<0$.

If a copy of T6 has a missing spindle in 2 directions, then it gives away to
its missing spindles at least $\frac28$.  Using Claim~\ref{cl5}(a) twice, we
see that the copy of T6 has either 4 trivial spindles on one side or else 2
trivial spindles (in the same direction) on each side.  In the first case, by
Claim \ref{cl5}(c) and (R5), the copy of T6 gives $\frac15$ to the copy of T1. 
In the second case, by Claim~\ref{cl5}(b) and (R5), the copy of T6 gives away
$\frac1{10} + \frac1{10}$.  Hence the copy of T6 finishes with excess at most
$\frac25-\frac14-\frac15 <0$.

If a copy of T6 has a missing spindle in 1 direction, then it gives away to
that missing spindle $\frac18$.  Also, the copy of T6 has 4 trivial spindles on
one side and 2 on the other. Hence, from (R5) it gives away at least $\frac15 +
\frac1{10}$ to adjacent copies of T1.  So, it finishes with excess weight at
most $\frac25-\frac18-\frac15 -\frac1{10}<0$.  
Finally, suppose the copy of T6 has no missing spindles.  Now by
Claim~\ref{cl5}(c), the copy of T6 is bordered on each side by a copy of T1
with only one T6 neighbor.  Hence, by (R5), the copy of T6 gives away
$\frac15+\frac15$, so finishes with excess at most 0.
\end{proof}
Claim~\ref{cl6} completes the proof that we outlined prior to the discharging phases.
Hence, $\chi_f(\Re^2)\ge\frac{76}{21}$.
\end{proof}

\begin{figure}
\begin{tikzpicture}[scale = 7]
\tikzstyle{VertexStyle} = []
\tikzstyle{EdgeStyle} = []
\tikzstyle{labeledStyle}=[shape = circle, minimum size = 6pt, inner sep = 1.2pt, draw]
\tikzstyle{unlabeledStyle}=[shape = circle, minimum size = 6pt, inner sep = 1.2pt, draw, fill]

\begin{scope}[shift={(-0.12,0)}]
\Vertex[style = _BlueShadeUpRightStyle, x = 1.12, y = 0.522, L = \tiny {}]{v44};
\Vertex[style = _BlueShadeDownLeftStyle, x = 1.05, y = 0.484, L = \tiny {}]{v44};

\Vertex[style = _BlueShadeUpRightStyle, x = 1.26, y = 0.524, L = \tiny {}]{v44};
\Vertex[style = _BlueShadeDownLeftStyle, x = 1.19, y = 0.484, L = \tiny {}]{v44};

\Vertex[style = _BlueShadeUpRightStyle, x = 1.19, y = 0.404, L = \tiny {}]{v44};
\Vertex[style = _BlueShadeDownLeftStyle, x = 1.12, y = 0.364, L = \tiny {}]{v44};

\Vertex[style = _BlueShadeUpRightStyle, x = 1.05, y = 0.404, L = \tiny {}]{v44};
\Vertex[style = _BlueShadeDownLeftStyle, x = 0.98, y = 0.364, L = \tiny {}]{v44};

\Vertex[style = _BlueShadeUpRightStyle, x = 1.33, y = 0.6453, L = \tiny {}]{v44};
\Vertex[style = _BlueShadeDownLeftStyle, x = 1.26, y = 0.604, L = \tiny {}]{v44};

\end{scope}

\Vertex[style = unlabeledStyle, x = 0.792115384615385, y = 0.807660256410256, L = \tiny {}]{v0}
\Vertex[style = unlabeledStyle, x = 0.932115384615384, y = 0.807660256410256, L = \tiny {}]{v1}
\Vertex[style = unlabeledStyle, x = 1.07211538461538, y = 0.807660256410256, L = \tiny {}]{v2}
\Vertex[style = unlabeledStyle, x = 1.21211538461538, y = 0.807660256410256, L = \tiny {}]{v3}
\Vertex[style = unlabeledStyle, x = 1.35211538461538, y = 0.807660256410256, L = \tiny {}]{v4}
\Vertex[style = unlabeledStyle, x = 0.722115384615384, y = 0.686360256410256, L = \tiny {}]{v5}
\Vertex[style = unlabeledStyle, x = 0.862115384615385, y = 0.686360256410256, L = \tiny {}]{v6}
\Vertex[style = unlabeledStyle, x = 1.00211538461538, y = 0.686360256410256, L = \tiny {}]{v7}
\Vertex[style = unlabeledStyle, x = 1.14211538461538, y = 0.686360256410256, L = \tiny {}]{v8}
\Vertex[style = unlabeledStyle, x = 1.28211538461538, y = 0.686360256410256, L = \tiny {}]{v9}
\Vertex[style = unlabeledStyle, x = 1.42211538461538, y = 0.686360256410256, L = \tiny {}]{v10}
\Vertex[style = unlabeledStyle, x = 0.652115384615385, y = 0.565160256410256, L = \tiny {}]{v11}
\Vertex[style = unlabeledStyle, x = 0.792115384615385, y = 0.565160256410256, L = \tiny {}]{v12}
\Vertex[style = unlabeledStyle, x = 0.932115384615384, y = 0.565160256410256, L = \tiny {}]{v13}
\Vertex[style = _RedDotStyle, x = 1.07211538461538, y = 0.565160256410256, L = \tiny {}]{v14}
\Vertex[style = unlabeledStyle, x = 1.21211538461538, y = 0.565160256410256, L = \tiny {}]{v15}
\Vertex[style = unlabeledStyle, x = 1.35211538461538, y = 0.565160256410256, L = \tiny {}]{v16}
\Vertex[style = unlabeledStyle, x = 1.49211538461538, y = 0.565160256410256, L = \tiny {}]{v17}
\Vertex[style = unlabeledStyle, x = 0.582115384615384, y = 0.443860256410256, L = \tiny {}]{v18}
\Vertex[style = unlabeledStyle, x = 0.722115384615384, y = 0.443860256410256, L = \tiny {}]{v19}
\Vertex[style = unlabeledStyle, x = 0.862115384615385, y = 0.443860256410256, L = \tiny {}]{v20}
\Vertex[style = unlabeledStyle, x = 1.00211538461538, y = 0.443860256410256, L = \tiny {}]{v21}
\Vertex[style = unlabeledStyle, x = 1.14211538461538, y = 0.443860256410256, L = \tiny {}]{v22}
\Vertex[style = unlabeledStyle, x = 1.28211538461538, y = 0.443860256410256, L = \tiny {}]{v23}
\Vertex[style = unlabeledStyle, x = 1.42211538461538, y = 0.443860256410256, L = \tiny {}]{v24}
\Vertex[style = unlabeledStyle, x = 1.56211538461538, y = 0.443860256410256, L = \tiny {}]{v25}
\Vertex[style = unlabeledStyle, x = 0.652115384615385, y = 0.322660256410256, L = \tiny {}]{v26}
\Vertex[style = _RedDotStyle, x = 0.792115384615385, y = 0.322660256410256, L = \tiny {}]{v27}
\Vertex[style = unlabeledStyle, x = 0.932115384615384, y = 0.322660256410256, L = \tiny {}]{v28}
\Vertex[style = unlabeledStyle, x = 1.07211538461538, y = 0.322660256410256, L = \tiny {}]{v29}
\Vertex[style = unlabeledStyle, x = 1.21211538461538, y = 0.322660256410256, L = \tiny {}]{v30}
\Vertex[style = unlabeledStyle, x = 1.35211538461538, y = 0.322660256410256, L = \tiny {}]{v31}
\Vertex[style = unlabeledStyle, x = 1.49211538461538, y = 0.322660256410256, L = \tiny {}]{v32}
\Vertex[style = unlabeledStyle, x = 0.722115384615384, y = 0.201360256410256, L = \tiny {}]{v33}
\Vertex[style = unlabeledStyle, x = 0.862115384615385, y = 0.201360256410256, L = \tiny {}]{v34}
\Vertex[style = _RedDotStyle, x = 1.00211538461538, y = 0.201360256410256, L = \tiny {}]{v35}
\Vertex[style = unlabeledStyle, x = 1.14211538461538, y = 0.201360256410256, L = \tiny {}]{v36}
\Vertex[style = unlabeledStyle, x = 1.28211538461538, y = 0.201360256410256, L = \tiny {}]{v37}
\Vertex[style = unlabeledStyle, x = 1.42211538461538, y = 0.201360256410256, L = \tiny {}]{v38}
\Vertex[style = unlabeledStyle, x = 0.792115384615385, y = 0.0801602564102564, L = \tiny {}]{v39}
\Vertex[style = unlabeledStyle, x = 0.932115384615384, y = 0.0801602564102564, L = \tiny {}]{v40}
\Vertex[style = unlabeledStyle, x = 1.07211538461538, y = 0.0801602564102564, L = \tiny {}]{v41}
\Vertex[style = unlabeledStyle, x = 1.21211538461538, y = 0.0801602564102564, L = \tiny {}]{v42}
\Vertex[style = unlabeledStyle, x = 1.35211538461538, y = 0.0801602564102564, L = \tiny {}]{v43}
\Edge[](v1)(v0)
\Edge[](v2)(v1)
\Edge[](v3)(v2)
\Edge[](v4)(v3)
\Edge[](v5)(v0)
\Edge[](v6)(v0)
\Edge[](v6)(v1)
\Edge[](v6)(v5)
\Edge[](v7)(v1)
\Edge[](v7)(v2)
\Edge[](v7)(v6)
\Edge[](v8)(v2)
\Edge[](v8)(v3)
\Edge[](v8)(v7)
\Edge[](v9)(v3)
\Edge[](v9)(v4)
\Edge[](v9)(v8)
\Edge[](v10)(v4)
\Edge[](v10)(v9)
\Edge[](v11)(v5)
\Edge[](v12)(v5)
\Edge[](v12)(v6)
\Edge[](v12)(v11)
\Edge[](v13)(v6)
\Edge[](v13)(v7)
\Edge[](v13)(v12)
\Edge[](v14)(v7)
\Edge[](v14)(v8)
\Edge[](v14)(v13)
\Edge[](v15)(v8)
\Edge[](v15)(v9)
\Edge[](v15)(v14)
\Edge[](v16)(v9)
\Edge[](v16)(v10)
\Edge[](v16)(v15)
\Edge[](v17)(v10)
\Edge[](v17)(v16)
\Edge[](v18)(v11)
\Edge[](v19)(v11)
\Edge[](v19)(v12)
\Edge[](v19)(v18)
\Edge[](v20)(v12)
\Edge[](v20)(v13)
\Edge[](v20)(v19)
\Edge[](v21)(v13)
\Edge[](v21)(v14)
\Edge[](v21)(v20)
\Edge[](v22)(v14)
\Edge[](v22)(v15)
\Edge[](v22)(v21)
\Edge[](v23)(v15)
\Edge[](v23)(v16)
\Edge[](v23)(v22)
\Edge[](v24)(v16)
\Edge[](v24)(v17)
\Edge[](v24)(v23)
\Edge[](v25)(v17)
\Edge[](v25)(v24)
\Edge[](v26)(v18)
\Edge[](v26)(v19)
\Edge[](v27)(v19)
\Edge[](v27)(v20)
\Edge[](v27)(v26)
\Edge[](v28)(v20)
\Edge[](v28)(v21)
\Edge[](v28)(v27)
\Edge[](v29)(v21)
\Edge[](v29)(v22)
\Edge[](v29)(v28)
\Edge[](v30)(v22)
\Edge[](v30)(v23)
\Edge[](v30)(v29)
\Edge[](v31)(v23)
\Edge[](v31)(v24)
\Edge[](v31)(v30)
\Edge[](v32)(v24)
\Edge[](v32)(v25)
\Edge[](v32)(v31)
\Edge[](v33)(v26)
\Edge[](v33)(v27)
\Edge[](v34)(v27)
\Edge[](v34)(v28)
\Edge[](v34)(v33)
\Edge[](v35)(v28)
\Edge[](v35)(v29)
\Edge[](v35)(v34)
\Edge[](v36)(v29)
\Edge[](v36)(v30)
\Edge[](v36)(v35)
\Edge[](v37)(v30)
\Edge[](v37)(v31)
\Edge[](v37)(v36)
\Edge[](v38)(v31)
\Edge[](v38)(v32)
\Edge[](v38)(v37)
\Edge[](v39)(v33)
\Edge[](v39)(v34)
\Edge[](v40)(v34)
\Edge[](v40)(v35)
\Edge[](v40)(v39)
\Edge[](v41)(v35)
\Edge[](v41)(v36)
\Edge[](v41)(v40)
\Edge[](v42)(v36)
\Edge[](v42)(v37)
\Edge[](v42)(v41)
\Edge[](v43)(v37)
\Edge[](v43)(v38)
\Edge[](v43)(v42)
\Edge[style = _GreenEdgeStyle](v14)(v27)
\Edge[style = _GreenEdgeStyle](v14)(v35)
\Edge[style = _GreenEdgeStyle](v35)(v27)

\begin{scope}[shift={(-.117,.005)}]
\draw[style = _ArrowStyle] (1.049,0.475) -- (1.14,0.533);
\draw[style = _ArrowStyle] (1.189,0.475) -- (1.28,0.533);
\draw[style = _ArrowStyle] (0.979,0.355) -- (1.07,0.413);
\draw[style = _ArrowStyle] (1.119,0.355) -- (1.21,0.413);
\draw[style = _ArrowStyle] (1.259,0.595) -- (1.35,0.653);
\end{scope}

\end{tikzpicture}
\caption{An additional 5 spindle block, used in the proof that
$\chi_f(\Re^2)\ge\frac{105}{29}$.\label{basetotip}
}
\end{figure}
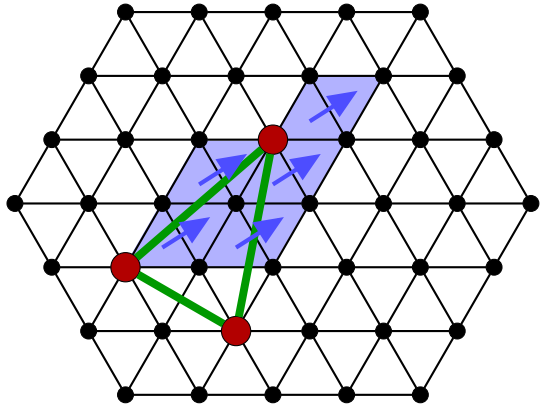

The lower bound $\frac{76}{21}\approx3.6190$ can be improved slightly to
$\frac{105}{29}\approx3.6207$, as follows.
We replace $\frac13$ in (R1) with $\frac{13}{42}$ and replace $\frac3{10}$ in
(R4) with $\frac27$.  A proof nearly identical to that above shows that every
tile (and spindle) finishes with excess at most 0, except for possibly copies of
T4.  In order to accomodate copies of T4, we need to consider another type of
5-spindle block, shown in Figure~\ref{basetotip}.  Although this proof can be
completed, the details are surprisingly complicated (the proof requires an
additional 4 phases), so we decided the reader would benefit more from the
version presented here.

%
%
%
%

The bound $\frac{105}{29}$ seems to be the best possible using this method.  The
key obstruction to further improvements is that we need weight at least 6 on each
core vertex.  If a core vertex has less weight than the weight on an independent
set among its spindle neighbors, then we cannot justify the assumption that $I$
intersects the core vertices in a maximal independent set.  This leads us to ask
for the value of $\lim_{d\to\infty}\chi_f(G'_d)$.  Specifically, is it larger
than $\frac{105}{29}$?

\section{Coloring Variants and Open Questions}

In this section, we discuss a few variants on the problem of coloring (or
fractionally coloring) the plane.  Typically, we either put restrictions on each
color class, or else we consider only a subset of the plane.  We also include a
few more open questions.

Falconer \cite{falconer1981} proved that if we require color classes to be
Lebesgue measurable, then any proper coloring of the plane requires at least
$5$ colors.  The existence of non-measurable sets depends on the axioms of set
theory we choose.  In particular, when we drop the Axiom of Choice, there are
(appealing) extensions of ZF in which all sets of reals are Lebesgue measurable
\cite{solovay1970}.  So, basically, Falconer's result says that we can
never hope to construct a $4$-coloring of the plane.  In the case of fractional
coloring, Sz{\'e}kely \cite{szekely1984} proved in 1984 that a fractional
coloring in which only measurable sets get non-zero weight must use total
weight at least $\frac{43}{12} \approx 3.5833$. Recently, this was improved by
Oliveira Filho and Vallentin \cite{de2010fourier} to $\approx 3.725$.

\begin{question}
Do the fractional chromatic number of the plane and the measurable fractional
chromatic number of the plane differ assuming ZFC?
\end{question}

The \emph{$j$-fold chromatic number} $\chi_j(G)$ of a graph $G$ is the least
integer $k$ so that it is possible to assign each vertex a set of $j$ elements
from $\{1,\ldots,k\}$ such that adjacent vertices get disjoint sets.  Clearly,
$\chi_f(G)\le \frac{\chi_j(G)}j$, since we can assign each of the $k$ color
classes weight $\frac1j$.  Recently, Grytczuk et.al.~\cite{GJSW14} studied
$j$-fold coloring of the plane.  Among other results, they generalized the
construction of Hochberg and O'Donnell to give good $j$-fold colorings for
small values of $j$.

Interesting results have been proved about chromatic numbers of extensions of
$\Q$, the first being Woodall's result \cite{woodall1973} that $\Q\times
\Q$ is 2-colorable.  Among other results, Fischer \cite{fischer1990} showed that
$\chi\parens{\Q(\sqrt{3}) \times \Q(\sqrt{3})} = 3$ and $\chi\parens{\Q(\sqrt{11})
\times \Q(\sqrt{11})} \le 4$.  The graphs $G'_d$ that we constructed have all
vertices in $\Q(\sqrt{3},\sqrt{11})\times \Q(\sqrt{3},\sqrt{11})$, so we have
actually shown that $\chi_f\parens{\Q(\sqrt{3},\sqrt{11})\times \Q(\sqrt{3},\sqrt{11})}
\ge \frac{76}{21}$. The natural extension of our construction is to attach spindles
at all rotations that are integer multiples of $\cos^{-1}(\frac56)$, and this
graph is still contained in $\Q(\sqrt{3},\sqrt{11})\times \Q(\sqrt{3},\sqrt{11})$. 
Does this graph have larger fractional chromatic number?  We do not know.

\begin{question}
What is the fractional chromatic number of $\Q(\sqrt{3},\sqrt{11})\times
\Q(\sqrt{3},\sqrt{11})$?  What about its chromatic number?
\end{question}

Another intriguing direction of work is unit distance graphs with higher girth.
Erd\H{o}s~\cite{erdos75} asked for which $k\ge 3$ there exist unit distance
graphs with chromatic number 4.  In his Ph.D. thesis~\cite{odonnell-thesis},
Paul O'Donnell showed that the answer is all $k\ge 3$.
Mohar~\cite{mohar01} extended this question to chromatic number 5 and 6.  
Much like the chromatic number (and fractional chromatic number) of the plane,
this problem remains open. 

\bibliographystyle{abbrvplain}
\bibliography{fractional-plane}

\end{document}